\RequirePackage{fix-cm}
\documentclass{svjour3x}                      
\smartqed  
\usepackage{graphicx}

\usepackage[margin=3cm]{geometry}
\usepackage{xspace}
\usepackage{amssymb,amsmath,makecell,epsfig,listings,bm}
\usepackage{caption}
\usepackage{float}
\usepackage{color}
\usepackage{xcolor}
\usepackage{subfig}

\usepackage{enumerate}
\usepackage{empheq}
\usepackage{mathtools}
\usepackage{pgf,tikz}
\usetikzlibrary{arrows}
\usetikzlibrary{shapes.misc}\usetikzlibrary{arrows.meta}
\tikzset{>={Latex[width=2mm,length=2mm]}}
\tikzset{cross/.style={cross out, draw=black, minimum size=2*(#1-\pgflinewidth), inner sep=0pt, outer sep=0pt},
	cross/.default={1pt}}
\graphicspath{{0_Figures/}}

\DeclareMathOperator{\sign}{sgn}
\DeclareMathOperator{\diam}{diam}

\DeclareMathOperator{\supp}{supp}

\DeclareSymbolFont{AMSb}{U}{msb}{m}{n}
\DeclareMathSymbol{\N}{\mathbin}{AMSb}{"4E}
\DeclareMathSymbol{\Z}{\mathbin}{AMSb}{"5A}
\DeclareMathSymbol{\R}{\mathbin}{AMSb}{"52}
\DeclareMathSymbol{\Q}{\mathbin}{AMSb}{"51}
\DeclareMathSymbol{\C}{\mathbin}{AMSb}{"43}

\newcommand{\bfb}{\mathbf{b}}

\newcommand{\bfa}{\mathbf{a}}
\newcommand{\bfx}{\mathbf{x}}
\newcommand{\bfv}{\mathbf{v}}

\newcommand{\calC}{\mathcal{C}}
\newcommand{\calW}{\mathcal{W}}
\newcommand{\calR}{\mathcal{R}}

\newcommand{\calM}{\mathcal{M}}
\newcommand{\calP}{\mathcal{P}}

\newcommand{\VHNA}{V_N}

\newcommand{\pdonetext}[2]{{\partial#1}/{\partial #2}}

\newcommand{\imag}{\mathrm{i}}
\newcommand{\calL}{\mathcal{L}}
\newcommand{\pmax}{p}
\newcommand{\colSet}{\calC_M}

\newcommand{\bfp}{\mathbf{p}}

\newcommand{\numScreens}{N_\Gamma}
\newcommand{\numQuad}{N_{\mathrm{Q}}}
\newcommand{\resub}[1]{{#1}}

\begin{document}
	
	\newcommand{\xx}{\mathbf{x}}
	\newcommand{\XX}{\mathbf{X}}
	\newcommand{\yy}{\mathbf{y}}
	\newcommand{\rr}{\mathbf{r}}
	\newcommand{\uu}{\mathbf{u}}
	\newcommand{\UU}{\mathbf{U}}
	\newcommand{\ff}{\mathbf{f}}
	\newcommand{\FF}{\mathbf{F}}
	\newcommand{\LL}{\mathbf{L}}
	\newcommand{\nn}{\mathbf{n}}
	\newcommand{\TT}{\mathbf{t}}
	\newcommand{\ee}{\mathbf{e}}
	\newcommand{\dd}{\mathbf{d}}
	\newcommand{\hatn}{\hat{\mathbf{n}}}
	\newcommand{\hatt}{\hat{\mathbf{t}}}
	\newcommand{\hatx}{\hat{\mathbf{x}}}
	\newcommand{\haty}{\hat{\mathbf{y}}}
	\newcommand{\hatz}{\hat{\mathbf{z}}}
	\newcommand{\half}{\frac{1}{2}}
	\newcommand{\ddfrac}[2]{\frac{{\rm d} #1}{{\rm d} #2}}
	\newcommand{\pfrac}[2]{\frac{\partial #1}{\partial #2}}
	
	\newcommand{\bx}{\mathbf{x}}
	\newcommand{\by}{\mathbf{y}}
	\newcommand{\bn}{\mathbf{n}}
	\newcommand{\bd}{\mathbf{d}}
	\newcommand{\bs}{\mathbf{s}}
	\newcommand{\cS}{\mathcal{S}}
	\newcommand{\re}{{\rm e}}
	\newcommand{\ri}{{\rm i}}
	\newcommand{\rd}{{\rm d}}
	\newcommand{\rf}{\eqref}
	\newcommand{\tGamma}{\widetilde{\Gamma}}
	\newcommand{\norm}[2]{\|#1\|_{#2}}
	\newcommand{\normt}[2]{\|#1\|_{#2}}
	\newcommand{\ord}[1]{\mathcal{O}\left(#1\right)}
	
	\DeclarePairedDelimiter{\ceil}{\lceil}{\rceil}

\title{Fast 	hybrid numerical-asymptotic boundary element methods for high frequency screen and aperture problems based on least-squares collocation}
	
%\titlerunning{Collocation HNA BEM for scattering by screens}
\titlerunning{Fast collocation HNA BEM for high frequency screen and aperture problems}

\author{A.\ Gibbs        \and
        D.\ P.\ Hewett        \and
        D.\ Huybrechs        \\\and
        E.\ Parolin
}

\institute{A.\ Gibbs and D.\ P.\ Hewett \at
              Department of Mathematics, University College London, London, UK \\
              Tel.: +44 20 7679 3935\\
              \email{andrew.gibbs@ucl.ac.uk}           
           \and
           D.\ Huybrechs  \at
            Department of Computer Science, KU Leuven, Leuven, Belgium
            \and
            E.\ Parolin \at
            Institut Polytechnique de Paris, Palaiseau, France
}


\maketitle

\begin{abstract}
	We present a hybrid numerical-asymptotic (HNA) boundary element method (BEM) for high frequency scattering by two-dimensional screens and apertures, whose computational cost to achieve any prescribed accuracy remains bounded with increasing frequency.
	Our method is a collocation implementation of the high order $hp$ HNA approximation space of Hewett et al. IMA J. Numer. Anal. 35 (2015), pp.\ 1698-1728, where a Galerkin implementation was studied. 
	An advantage of the current collocation scheme is that the one-dimensional highly oscillatory singular integrals appearing in the BEM matrix entries are significantly easier to evaluate than the two-dimensional integrals appearing in the Galerkin case, which leads to much faster computation times. 
	Here we compute the required integrals at frequency-independent cost using the numerical method of steepest descent, which involves complex contour deformation. 
	The change from Galerkin to collocation is nontrivial because naive collocation implementations based on square linear systems suffer from severe numerical instabilities associated with the numerical redundancy of the HNA basis, which produces highly ill-conditioned BEM matrices. In this paper we show how these instabilities can be removed by oversampling, and solving the resulting overdetermined collocation system in a weighted least-squares sense using a truncated singular value decomposition. On the basis of our numerical experiments, the amount of oversampling required to stabilise the method is modest (around 25\% typically suffices), and independent of frequency. 
	As an application of our method we present numerical results for high frequency scattering by prefractal approximations to the middle-third Cantor set.
\end{abstract}

\section{Introduction}

The numerical simulation of high frequency (short wavelength) acoustic and electromagnetic scattering is challenging because of the need to capture the rapid oscillations in the wave field. Conventional Boundary Element Methods (BEMs), based on piecewise-polynomial basis functions, are computationally expensive because they require a fixed number of degrees of freedom (DOFs) per wavelength to achieve accurate solutions. This leads to very large (dense) BEM matrices which are costly to store and invert. 

By contrast, \textit{hybrid numerical-asymptotic} (HNA) BEMs 
use basis functions built from piecewise polynomials on coarse meshes multiplied by certain oscillatory functions, 
chosen based on partial knowledge of the high frequency asymptotic solution behaviour, as described by Geometrical Optics (GO) and the Geometrical Theory of Diffraction (GTD) \cite{Keller,Kinber,james1986geometrical}. 
The goal of the HNA approach (reviewed in \cite{acta}) is to achieve a fixed accuracy of approximation using a number of DOFs that is relatively small and frequency-independent, or only modestly (e.g.\ logarithmically) frequency-dependent, making it easier to store and invert the BEM matrix at high frequencies.
HNA BEMs have been successfully developed for  
scattering by 
impenetrable convex \cite{CWL:07,hewett2013high,CWLM}, nonconvex \cite{chandler2012high,hewett:shadow} and curvilinear \cite{LaMoCh:10} polygons, two-dimensional screens and apertures \cite{hewett2014frequency}, smooth convex two-dimensional \cite{dominguez2007hybrid,bruno2004prescribed,AsHu:14,EcOz:17,EcUr:16} and three-dimensional \cite{GH11} scatterers, three-dimensional rectangular screens \cite{Hargreaves}, penetrable convex polygons \cite{groth2015hybrid,groth2018hybrid} and, recently, for certain multi-obstacle scattering problems \cite{ChGiLaMo:19}. 

The use of oscillatory bases in HNA methods leads to an essentially frequency-independent BEM system size. However, the BEM matrix entries now involve highly oscillatory singular integrals, which need to be evaluated efficiently if one is to achieve the ``holy grail'' of frequency-independent computational cost. 
The majority of HNA methods in the literature to date 
(e.g.\ \cite{CWL:07,hewett2013high,CWLM,chandler2012high,hewett:shadow,LaMoCh:10,hewett2014frequency,dominguez2007hybrid,GH11,Hargreaves,EcOz:17,EcUr:16,groth2018hybrid}) 
are based on Galerkin discretisations. This leads to provably stable approximations and, in many cases (e.g.\ \cite{CWL:07,hewett2013high,chandler2012high,hewett:shadow,hewett2014frequency,EcOz:17,EcUr:16}), provides the framework for a rigorous, frequency-explicit convergence analysis. 
But Galerkin testing produces high-dimensional oscillatory integrals, which complicates the development of fast quadrature routines. 

In order to develop HNA methods for more practically relevant problems than those tackled so far (in particular, three-dimensional problems), 
it would be highly advantageous to be able to use a collocation or Nystr\"om approach, for which numerical quadrature is simpler. 
Existing work in this direction includes for example the piecewise-constant collocation method for convex polygons in \cite{arden2007collocation}, the B-spline method for two-dimensional smooth convex scatterers in \cite{Gi:07} (which is based on the earlier technical report \cite{giladi2004asymptotically}) and the closely related Nystr\"om method in \cite{bruno2004prescribed}. 
However, these approaches are not generally supported by rigorous stability and convergence analysis, a practical consequence of which being that it is not at all obvious how to determine collocation/quadrature point distributions leading to stable approximations. This question is particularly delicate when working with non-smooth scatterers and HNA approximation spaces built with overlapping meshes to represent different high frequency solution components (as in \cite{arden2007collocation}), which can exhibit a high degree of numerical ``redundancy'' and hence severe ill-conditioning (see e.g. the discussion in \cite[\S2]{arden2007collocation} and \S\ref{sec:SVD} below).

Our aim in this paper is to demonstrate that accurate and efficient collocation-based HNA methods can indeed be developed for high frequency scattering problems involving non-smooth scatterers. The key new idea, apparently not employed previously in HNA methods, is to stabilise the method using oversampling, and solve the resulting overdetermined collocation linear system in a weighted least-squares sense.

\begin{figure}
	\centering
	\includegraphics[width=\linewidth]{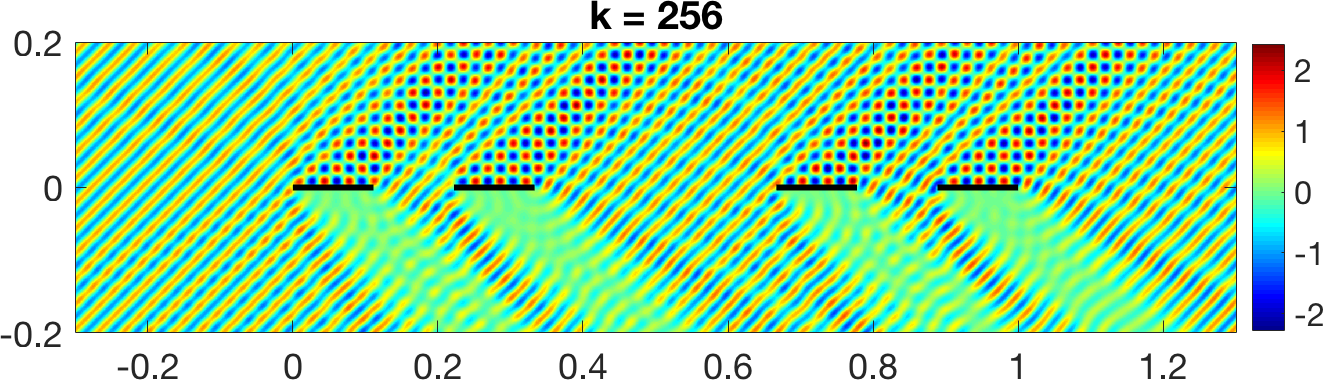}
	\caption{Real part of the total field for scattering by 
		$N_\Gamma=4$	
		%	four 
		colinear sound-soft screens constituting the order $j=2$ prefractal of the middle-third Cantor set (for more details see \S\ref{sec:Cantor}). The incident wave propagates from top left to bottom right. 
	}
	\label{fig:mutliscreenplot}
\end{figure}

We present our oversampled collocation approach in the context of a specific model problem, namely 
high frequency two-dimensional scattering by colinear screens and apertures.
For an illustration of the problem see Fig.~\ref{fig:mutliscreenplot}. 
%,as studied recently in \cite{priddin2019applying}, where an iterative Wiener-Hopf approaches was considered.  
Such problems have numerous applications in acoustics, electromagnetics and water waves - for details we refer to the extensive reference list in \cite{hewett2014frequency}, noting also the recent work on iterative Wiener-Hopf approaches to this problem in \cite{priddin2019applying}.  
Our collocation HNA BEM for the screen problem uses the same high order HNA approximation space as the Galerkin method presented in \cite{hewett2014frequency}, which is based on an $hp$ refinement strategy, giving exponential accuracy with increasing polynomial degree (see Theorem \ref{dudnThm} below). 
To stabilise our collocation method we take more collocation points than degrees of freedom, and solve the resulting rectangular (overdetermined) system in a weighted least-squares sense, using a truncated singular value decomposition. The one-dimensional oscillatory singular integrals appearing in the BEM matrix are evaluated accurately and efficiently using the numerical method of steepest descent \cite{DeHu:09,HuVa:06,DaanBook}, which is based on complex contour deformation, combined with generalised Gaussian quadrature \cite{HuCo:08}, to handle the singularities without the need for mesh grading. 

By a series of numerical experiments we demonstrate that our HNA collocation BEM can achieve a fixed solution accuracy in a computation time that is bounded (in fact, sometimes even decreasing) as the frequency tends to infinity. A key empirical observation (that we are currently unable to explain theoretically) is that the \emph{oversampling threshold}, which governs the number of collocation points per DOF, does not need to increase with increasing frequency in order to maintain accuracy. 

This paper has as its genesis the MSc dissertation \cite{emile} of the fourth author, in which a range of HNA approximation spaces and collocation point distribution strategies were compared for square collocation BEM matrices. Numerical instabilities encountered in \cite{emile} motivated the investigation of the oversampled collocation method presented here. 
The oscillatory integration in \cite{emile} was carried out using Filon quadrature (see e.g. \cite{Iserles2004,Huybrechs2009,DaanBook}), combined with mesh grading to handle singularities (as in \cite{DoGrKi:13}). 

\resub{We point out that our oversampled collocation approach differs slightly from that of the CHIEF method of \cite{Sc:68}. The CHIEF method was developed for scattering by closed surfaces, and introduces extra collocation points in the \emph{interior} of the scatterer to stabilise integral equation formulations that are ill-posed at certain values of $k$, corresponding to interior resonances. In the current paper we are introducing extra collocation points on the scatterer \emph{boundary}, to stabilise an integral equation formulation for scattering by an open surface (in particular, the scatterer has no interior), which is well-posed for all $k$ but suffers from ill-conditioning when approximated numerically.}

The structure of the paper is as follows. 
In section \S\ref{sec:BVP} we define the scattering problem and its boundary integral equation reformulation. In \S\ref{sec:HNAspace} we describe the HNA approximation space of \cite{hewett2014frequency}, and recall the corresponding best approximation error estimate. In \S\ref{sec:collocation} we present our regularised collocation approach, and in \S\ref{sec:SVD} we discuss the philosophy behind the truncated SVD approach used to solve the least-squares problem. In \S\ref{sec:oscillatoryquadrature} we outline our numerical quadrature scheme for evaluating the oscillatory and singular integrals appearing in the BEM system. In \S\ref{sec:numResults} we present a series of numerical experiments demonstrating the performance of our method, including an application to scattering by the middle-third Cantor set. 

\section{Problem statement}\label{sec:BVP}

We consider 2D time-harmonic acoustic scattering by a sound-soft (Dirichlet) screen %$\Gamma$, 
\begin{align}
	\label{eqn:GammaDef}
	\Gamma=\{ (x_1,0)\in\R^2:x_1\in\tGamma\},
	\qquad
	\tGamma = \bigcup_{j=1}^{\numScreens} (s_{2j-1},s_{2j}),
\end{align}
a bounded and relatively open subset of $\Gamma_\infty:=\{ \bx=(x_1,x_2)\in\R^2:x_2=0\}$ comprising $\numScreens\geq 1$ colinear disjoint open intervals, with $0=s_1<s_2<\ldots<s_{2\numScreens}=\diam{\Gamma}$. 
We assume that lengths have been nondimensionalised so that $s_{2\numScreens}=\diam{\Gamma}=1$. 
For each $j=1,\ldots,\numScreens$ we set $\Gamma_j:=(s_{2j-1},s_{2j})\times\{0\}\subset\R^2$.  
The propagation domain is $D:=\R^2\setminus \overline\Gamma$,  
and on $\Gamma$ we define the normal vector $\bn=(0,1)$.

As the incident wave we take $u^i(\bx) :=\re^{\ri  k \bx\cdot \bd}$, 
where $k>0$ is the nondimensional wavenumber 
and $\bd=(d_1,d_2)\in\mathbb{R}^2$ is a unit direction vector. 
The boundary value problem (BVP) to be solved is: Find $u\in C^2\left(D\right)\cap C(\R^2)\cap  W^1_{\mathrm{loc}}(D)$ s.t.\ 
\begin{align}
	\Delta u+k^2u  & =  0, \qquad \mbox{in }D, \label{eqn:HE1} \\
	u & =  0, \qquad \mbox{on }\Gamma,\label{eqn:bc1}
\end{align}
with the scattered field $u^s:=u-u^i$ satisfying the Sommerfeld radiation condition (see, e.g., \cite[Equation (2.9)]{acta}). 
Here $W^1_{\rm loc}(D)$ is the usual space of locally square-integrable functions whose weak gradient is also locally square-integrable. Well-posedness of the BVP \rf{eqn:HE1}-\rf{eqn:bc1} is proved, for example, in \cite{StWe:84}. For an example solution see Fig.~\ref{fig:mutliscreenplot}.

\begin{remark}\label{rem:Babz}
	By Babinet's principle, the Dirichlet screen problem \rf{eqn:HE1}-\rf{eqn:bc1} is equivalent to its complementary aperture problem, in which $u^i$ impinges on an aperture $\Gamma$ in a sound-hard (Neumann) screen occupying $\Gamma_\infty\setminus\overline{\Gamma}$. (For a precise definition of the aperture BVP see \cite[Definition 1.2]{hewett2014frequency}.)
	Explicitly, with $U^\pm = \{\bx\in\R^2:\pm x_2>0\}$, $u^r(\bx) :=\re^{\ri  k \bx\cdot \bd'}$, $d_2<0$, and $\bd':=(d_1,-d_2)$, 
	the screen and aperture solutions $u$ and $u'$ are connected by the formula \cite[Equation (3.8)]{hewett2014frequency}
	\begin{align}
		\label{eqn:corres}
		u'(\bx )=
		\begin{cases}
			u^r(\bx) + u(\bx), & \bx\in U^+,\\
			u^i(\bx) - u(\bx), & \bx\in U^-.
		\end{cases}
	\end{align}
\end{remark}

Theorem \ref{DirEquivThm} below reformulates the BVP \rf{eqn:HE1}-\rf{eqn:bc1} as an an integral equation for the normal derivative jump $[\pdonetext{u}{\bn}]:=(\pdonetext{u}{\bn})^+ - (\pdonetext{u}{\bn})^-$, where $^\pm$ respectively denote the values on the top ($+$) and bottom ($-$) of $\Gamma$. 
We first clarify some notation (for more detail see \cite{hewett2014frequency}).
For $s\in\R$, we denote by $H^s(\R)$ the usual Sobolev space on $\R$, which, following \cite{hewett2014frequency,CoercScreen2}, we equip with the wavenumber-dependent norm
\begin{align}\label{def:kDepSobNorm}
	\norm{u}{H_k^{s}(\R)}^2 :=\int_{-\infty}^\infty(k^2+\xi^2)^{s}\,|\hat{u}(\xi)|^2\,\rd \xi,
\end{align}
where $\hat{u}$ denotes the Fourier transform of $u$. 
We set $\widetilde H^s(\widetilde\Gamma):=\overline{C^\infty_0(\widetilde\Gamma)}^{H^s(\R)}$, equipped with the norm $\norm{u}{H_k^{s}(\R)}$ inherited from $H^s(\R)$, and $H^s(\widetilde\Gamma):=\{u=U|_{\widetilde\Gamma}:U\in H^s(\R)\}$, equipped with the norm $\|u\|_{H^s_k(\widetilde\Gamma)} = \inf \|U\|_{H^s_k(\R)}$, %\agnote{should we specify what this infimum is over?}, 
where the infimum is taken over all $U\in H^s(\R)$ such that $u=U|_{\widetilde\Gamma}$.
We identify $\Gamma\subset\Gamma_\infty$ with $\widetilde\Gamma\subset\R$ in the natural way and define $H^s(\Gamma_\infty):=H^s(\R)$, $\widetilde H^s(\Gamma):=\widetilde H^s(\widetilde \Gamma)$, $H^s(\Gamma):=H^s(\widetilde \Gamma)$ etc.
We denote by $\cS:\widetilde{H}^{-1/2}(\Gamma)\to C^2(D)\cap W^1_{loc}(\R^2)$ and $S:\widetilde{H}^{-1/2}(\Gamma)\to H^{1/2}(\Gamma)$ the standard single-layer potential and single-layer boundary integral operator, which for $\phi\in L^2(\Gamma)$ have the integral representations
\begin{align*}
	\cS\phi (\bx) &= \int_\Gamma \Phi(\bx,\by) \phi(\by)\, \rd s(\by), \qquad \bx\in \R^2,\\
	S\phi (\bx) &= \int_\Gamma \Phi(\bx,\by) \phi(\by)\, \rd s(\by), \qquad \bx\in \Gamma,
\end{align*}
where $\Phi(\bx,\by) := (\ri/4)H_0^{(1)}(k|\bx-\by|)$ is the fundamental solution of \rf{eqn:HE1}.

\begin{theorem}[{\cite[Theorem 1.7]{StWe:84}, \cite[Theorem 3.1]{hewett2014frequency}}]
	\label{DirEquivThm}
	Suppose that $u$ is a solution of the BVP \rf{eqn:HE1}-\rf{eqn:bc1}. Then the representation formula 
	\begin{align}
		\label{eqn:D_RepThm}
		u(\bx )=u^i(\bx)  -\cS v(\bx), \qquad\bx\in D,
	\end{align}
	holds, with $v=[\pdonetext{u}{\bn}]\in \widetilde{H}^{-1/2}(\Gamma)$. Furthermore, $v$ satisfies the integral equation
	\begin{align}
		\label{BIE_sl}
		S v=f,
	\end{align}
	with $f:=u^i|_\Gamma\in H^{1/2}(\Gamma)$. 
	Conversely, suppose that $v\in \widetilde{H}^{-1/2}(\Gamma)$ satisfies \rf{BIE_sl}. 
	Then $u$ defined by \rf{eqn:D_RepThm} satisfies the BVP \rf{eqn:HE1}-\rf{eqn:bc1}, and $[\pdonetext{u}{\bn}]=v$.
\end{theorem}

\resub{We note that \eqref{BIE_sl} is well-posed for all $k>0$, with no spurious frequencies. Indeed, the operator $S$ is coercive (strongly elliptic) on $\widetilde H^{-1/2}(\Gamma)$ \cite{CoercScreen2}.}

\section{Hybrid numerical-asymptotic approximation space}\label{sec:HNAspace}

Our $hp$ HNA BEM approximation space for solving the BIE \eqref{BIE_sl} is identical to that considered in \cite{hewett2014frequency}. Its design is motivated by the following regularity result.

\begin{theorem}[{\cite[Theorem 4.1, Lemma 4.5]{hewett2014frequency}}]
	\label{vpmThm}
	The{ solution $v=[\pdonetext{u}{\bn}]$ of the BIE \eqref{BIE_sl}} can be decomposed on each screen component $\Gamma_j$, $j=1,\ldots,\numScreens$, as
	\begin{align}\label{Decomp}
		{v(\bx(s))}
		= \Psi(\bfx(s))+ v_j^+(s-s_{2j-1})\re^{\ri ks} +v_j^-(s_{2j}-s) \re^{-\ri ks}, & \quad s\in(s_{2j-1},s_{2j}), 
	\end{align}
	where $\Psi:=-2\sign(d_2)\pdonetext{u^i}{\bn}$,  
	and the functions $v_j^\pm(s)$ are analytic in $\Re{s}>0$, and, for each $k_0>0$, there exists a constant $C>0$, depending only on $k_0$, such that
	\begin{align}
		\label{vpmBounds}
		|v_j^\pm(s)|\leq C (1+k)k|ks|^{-1/2}, \qquad \Re{s}>0, \quad k\geq k_0.
	\end{align}
\end{theorem}

The known term $\Psi$ constitutes the ``physical optics'' approximation, describing the contribution to $v=[\pdonetext{u}{\bn}]$ of the incident and specularly reflected waves. 
The bounds \rf{vpmBounds} imply that the unknown functions $v_j^\pm$, which correspond to the amplitudes of the diffracted waves, are non-oscillatory as a function of $s$ (see \cite[Remark 4.2]{hewett2014frequency}), which means they can be approximated much more efficiently than the full (oscillatory) solution $v$. 
Hence, rather than approximating 
$v$ 
itself using piecewise polynomials on a fine ($k$-dependent) mesh (as in conventional BEMs), our HNA approximation space uses the decomposition~\rf{Decomp}, with $\Psi$ evaluated analytically and 
$v_j^+$ and $v_j^-$ replaced by piecewise polynomials on coarse ($k$-independent) meshes, graded to account for 
the singularities of $v_j^+$ at $s_{2j-1}$, and of $v_j^-$ at $s_{2j}$.

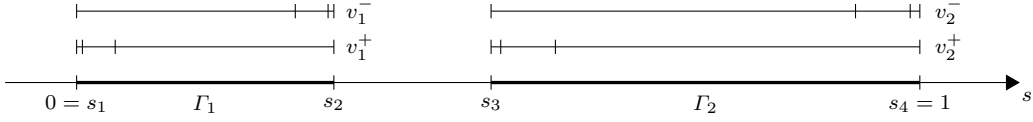
\begin{figure}[t!]
	\centering
	\begin{tikzpicture}[scale=0.95]
	\draw[line width=1.25pt] (3,0)--(-3,0);
	\draw (3,0.5)--(-3,0.5);
	\draw (3,1)--(-3,1);
	\foreach \x in {-3,-2.865,-2.1,3}
	{
		\draw (\x ,0.4) -- (\x ,0.6);
		\draw (-\x ,0.9) --(-\x,1.1);
	}
	\path (3.1,1)node[anchor=west]{$v_2^-$};
	\path (3.1,0.5)node[anchor=west]{$v_2^+$};
	\path (4.9,0.5)node[anchor=west]{\color{white}.};
	\path (0,-0.3)node{$\Gamma_2$};
	\begin{scope}[yshift=0.7cm]
	\draw (-3 ,-0.6) -- (-3 ,-0.8);
	\draw (3 ,-0.6) -- (3 ,-0.8);
	\path (-3,-1.0)node{$s_3$};
	\path (3,-1.0)node{$s_4=1$};
	\end{scope}
	\begin{scope}[x=0.6cm,y=1.0cm,xshift=-7cm]
	\draw[line width=1.25pt] (3,0)--(-3,0);
	\draw (3,0.5)--(-3,0.5);
	\draw (3,1)--(-3,1);
	\foreach \x in {-3,-2.865,-2.1,3}
	{
		\draw (\x ,0.4) -- (\x ,0.6);
		\draw (-\x ,0.9) --(-\x,1.1);
	}
	\path (3.1,1)node[anchor=west]{$v_1^-$};
	\path (3.1,0.5)node[anchor=west]{$v_1^+$};
	\path (0,-0.3)node{$\Gamma_1$};
	\begin{scope}[yshift=0.7cm]
	\draw (-3 ,-0.6) -- (-3 ,-0.8);
	\draw (3 ,-0.6) -- (3 ,-0.8);
	\path (-3,-1.0)node{$0=s_1$};
	\path (3,-1.0)node{$s_2$};
	\end{scope}
	\end{scope}
	\begin{scope}[yshift=0.7cm]
	\draw[arrows={-triangle 60}] (-9.8,-0.7)--(4.4,-0.7);
	\path (4.5,-0.9)node{$s$};
	\end{scope}
	\end{tikzpicture}
	\caption{Illustration of the overlapping graded meshes used to approximate the amplitudes $v_j^\pm$ in \rf{Decomp}, in the case where $\Gamma$ comprises two components, $\Gamma_1$ and $\Gamma_2$. %
	}
	\label{fig:grading}
\end{figure}

To describe the meshes we use, we first define a geometrically graded mesh \(\calM\) on the interval $[0,1]$ with $l$ %layers 
elements  
and grading parameter $\sigma\in(0,1/2)$, 
and an associated space of piecewise polynomials \(\calP_{\mathbf{p}}(\calM)\), 
by 
\begin{gather*}
	\calM:=\{x_i\}_{i=0}^l, \qquad x_0=0,\,\,x_i=\sigma^{l-i}, \,i=1,\ldots,l,\\
	{\calP}_{\mathbf{p}}(\calM):=\left\{ \rho:[0,1]\rightarrow\mathbb{C}:{\rho}|_{\left( x_{i-1}, x_i\right)} \text{ is a polynomial of degree $\leq$ } ({\mathbf{p}})_i,\, i=1,\ldots,l\right\},
\end{gather*}
where the degree vector $\bfp\in\N_0^{l}$ is defined for a given maximum degree $\pmax>0$ by
\begin{equation}\label{piDef1}
	({\bf{p}})_i:=
	\begin{cases}
		\pmax - \left\lfloor \frac{l+1-i}{l}\pmax \right\rfloor,& 1\leq i\leq l-1,\\
		\pmax, & i=l.\\
	\end{cases}
\end{equation}
Then for each screen segment $\Gamma_j$ of length $L_j:=s_{2j}-s_{2j-1}$, we define the meshes
\[
\calM^-_{j} := s_{2j-1} + L_j\calM,\qquad \calM^+_{j} := s_{2j} - L_j\calM,
\]
and denote by ${\calP}_{\mathbf{p}}(\calM_j^\pm)$ the associated spaces of piecewise polynomials, defined analogously to ${\calP}_{\mathbf{p}}(\calM)$ above. For an illustration of the resulting meshes see Fig.~\ref{fig:grading}.
Finally, our HNA approximation space for approximating the difference 
$v-\Psi$ 
is 
\begin{equation}\label{HNASpace}
	\VHNA := \bigoplus_{j=1}^{\numScreens}\left[\calP_{\mathbf{p}}(\calM_j^+)\re^{\ri k s} +  \calP_{\mathbf{p}}(\calM_j^-)\re^{-\ri k s}\right],
\end{equation}
with dimension 
\[
N = 2\numScreens\sum_{i=1}^{l}[(\bfp)_i+1].
\]

We then have the following best approximation error estimate. 
\begin{theorem}[{\cite[Theorem 5.1]{hewett2014frequency}}]
	\label{dudnThm}
	Let $k_0>0$ and $l\geq cp$ for some constant $c>0$. 
	Then for every $\delta>0$ there exists a constant $C>0$, depending only on $\delta$, $\sigma$, $\numScreens$ and $c$, and a constant $\tau>0$, depending only on $\delta$, $\sigma$ and $c$, such that
	\begin{align}
		\label{BestAppdudn}
		\inf_{\resub{\eta}\in \VHNA}\norm{(v-\Psi)-\resub{\eta}}{\widetilde{H}_k^{-1/2}(\Gamma)}
		\leq C
		(1+k)
		k^\delta\,\re^{-\pmax\tau}, \qquad k\geq k_0.
	\end{align}
\end{theorem}
Theorem \ref{dudnThm} states that $\VHNA$ can approximate 
$v - \Psi$ 
with an error which decreases exponentially as the maximum polynomial degree $\pmax$ tends to infinity, with a $k$-independent exponent. 
\resub{While the constant premultiplying the exponential factor grows with increasing $k$, it does so only modestly (like $O(k^{1+\delta}$))}, meaning that we can achieve any specified accuracy of approximation with $\pmax$ growing only logarithmically in $k$ as $k\to\infty$. This corresponds to the number of degrees of freedom $N$ growing like $\log^2{k}$ as $k\to\infty$. 
We shall see later that this logarithmic growth appears unnecessary in practice.

\section{Collocation method}\label{sec:collocation}

The difference 
\[ 
\phi:= v - \Psi
\]
satisfies the BIE
\begin{equation}\label{eq:HNABIE}
	S\phi = f - S\Psi.
\end{equation}
Let $\{\varphi_n\}_{n=1}^N$ denote a basis for $\VHNA$. 
In our implementation, each basis function $\varphi_n$ consists of an exponential $\re^{\ri ks}$ or $\re^{- \ri ks}$ multiplied by an appropriately scaled and translated Legendre polynomial
supported on a single element of one of the meshes $\calM_j^+$ or $\calM_j^-$ for some $j\in\{1,\ldots,N_\Gamma\}$, normalised so that $\|\varphi_n\|_{L^2(\Gamma)}=1$. 

To select an element $\sum_{n=1}^N a_n\varphi_n \in \VHNA$ approximating $\phi \in \tilde{H}^{-1/2}(\Gamma)$ we 
ask that the BIE \rf{eq:HNABIE}, with $\phi$ replaced by $\sum_{n=1}^N a_n\varphi_n \in \VHNA$, should hold at a set of $M$ collocation points 
\[
\colSet:=\{c_1,\ldots,c_M\}\subset \widetilde{\Gamma}=\bigcup_{j=1}^{\numScreens}(s_{2j-1},s_{2j}).
\]
As we shall see, for reasons of stability it is useful to oversample, taking $M> N$. 
The resulting system of linear equations for the coefficients $a_n$, $n=1,\ldots,N$, is then overdetermined, so we seek a weighted least-squares solution. %In more detail, 
Explicitly, given a set of weights $\calW:=\{w_m\}_{m=1}^M$ corresponding to the collocation points $\colSet$ (our choice of weights will be detailed below), we define $A\in \C^{M\times N}$ and $\bfb\in \C^M$ by
\begin{align}
	A_{m,n} := \sqrt{w_m}S\varphi_n(c_m),\quad b_m := \sqrt{w_m}\left(f(c_m)-S\Psi(c_m)\right), \quad &m=1,\ldots,M, \notag\\& n=1,\ldots, N,
	\label{eq:systemDeets}
\end{align}
and then seek $\bfa=(a_1,\ldots,a_n)^T\in\C^N$ 
minimising the residual
\begin{equation}\label{eq:LS}
	\calR:=\|A\bfa - \bfb \|_{2}.
\end{equation}
This problem is highly ill-conditioned, but can be regularised using a truncated singular value decomposition (SVD), as in e.g.\ \cite{engl1996regularization,hansen2013least,neumaier1998solving}. 
The full SVD of $A$ takes the form
\[
A= U\Sigma V^*,
\]
where $U$ and $V$ are unitary matrices, $V^*$ denotes the conjugate transpose of $V$, and $\Sigma $ is an $M\times N$ diagonal matrix.
Denote by $\sigma_n$ the $n$th entry of the diagonal of $\Sigma$. 
To regularise, we introduce a small threshold $\epsilon>0$, and define a modified diagonal matrix $\Sigma^{\epsilon}$, which has $n$th entry equal to $\sigma_n$ if $\sigma_n/ \max_{n'=1,\ldots,N}\sigma_{n'}>\epsilon $, and zero otherwise. This provides a regularisation 
\begin{equation}\label{eq:regSVD}
	A^{\epsilon}  = U\Sigma^{\epsilon} V^*,
\end{equation}
from which we can define $\bfa$ via a pseudo-inverse as
\begin{equation}\label{eq:solution_pinv}
	\bfa = (A^\epsilon)^\dag \bfb= V(\Sigma^\epsilon)^\dag U^*\bfb,
\end{equation}
where $(\Sigma^\epsilon)^\dag$ is the (diagonal) pseudo-inverse of $\Sigma^{\epsilon}$ with entries $\sigma^\dag_n$, defined such that $\sigma^\dag_n=1/\sigma_n$ if $\sigma_n/\max_{n'=1,\ldots,N}\sigma_{n'}>\epsilon $, and $\sigma^\dag_n=0$ otherwise.

Regarding the placement of collocation points, while we have no theoretical results to guide us, as a result of extensive numerical experiments (including those carried out for square systems ($M=N$) in \cite{emile}), we arrived at the following prescription.

We first select an oversampling threshold $C_\text{OS}>1$. 
For each $j\in \{1,\ldots,N_\Gamma\}$ let $s_{2j-1}=x^\pm_{j,0}<x^\pm_{j,1}<\cdots < x^\pm_{j,l} = s_{2j}$ denote the points of the overlapping meshes $\calM_j^\pm$ on $\Gamma_j$ (cf.\ Fig.~\ref{fig:grading}). On each element $[x^\pm_{j,i-1},x^\pm_{j,i}]$ we allocate $M^\pm_{j,i}:={\lceil{C_\text{OS}(p^\pm_{j,i} + 1)}\rceil}$ collocation points, where $p^\pm_{j,i}$ is the largest polynomial degree included on $\VHNA$ on that element. For the elements 
$[x^+_{j,0},x^+_{j,1}]$, \ldots $[x^+_{j,l-2},x^+_{j,l-1}]$ 
and
$[x^-_{j,1},x^-_{j,2}]$, \ldots $[x^-_{j,l-1},x^-_{j,l}]$ 
we place the collocation points at the first kind Chebyshev nodes within that element.
If we were to follow this same prescription for the largest elements, $[x^+_{j,l-1},x^+_{j,l}]$ and $[x^-_{j,0},x^-_{j,1}]$, 
the overlapping nature of the meshes would mean that collocation points could coincide with, or lie very close to other collocation points already selected on the smaller elements. To avoid this, for these largest elements we allocate $M^+_{j,1}+M^-_{j,l}=2{\lceil{C_\text{OS}(p + 1)}\rceil}$ collocation points at the first kind Chebyshev nodes in the interval $[x^+_{j,l-1},x^-_{j,1}]$ (the intersection of the two largest elements; recall that $0<\sigma<1/2$ so this intersection is non-trivial). The total number of collocation points is then
\begin{equation}\label{eq:MdivN_ne_Cos}
	M = \sum_{j=1}^{N_\Gamma}\sum_{\pm}\sum_{i=1}^{l} M^\pm_{j,i} 
	=\sum_{j=1}^{N_\Gamma}\sum_{\pm}\sum_{i=1}^{l} \lceil{C_\text{OS}(p^\pm_{j,i} + 1)}\rceil \approx C_\text{OS} N, \end{equation}
and our prescription ensures that
%Note that if $C_\text{OS}$ is an integer then $M=C_\text{OS}N$, but if $C_\text{OS}$ is not an integer then $M>C_\text{OS}N$ in general. 
\begin{align}
	\label{eq:MoverN}
	M \geq C_\text{OS}N, \quad \text{with equality if } C_\text{OS} \text{ is an integer.}
\end{align}

In \rf{eq:systemDeets} we choose weights $\{w_m\}_{m=1}^M$ that are inversely proportional to the approximate local density of collocation points. %, to compensate for the clustering of collocation points near the screen endpoints (due to the geometric mesh grading). 
Explicitly, suppose that $[a,b]$ is a mesh element on which we have allocated $M_{[a,b]}$ collocation points at the first-kind Chebyshev nodes
\[
c_\ell  = (a+b)/2 + (b-a)\cos(\pi(\ell -1/2)/M_{[a,b]})/2,\quad \ell =1,\ldots,M_{[a,b]}.
\]
Then to the points $\{c_\ell \}$ we assign the weights $\{w_\ell \}$ defined by
\begin{align}
	&w_\ell  = (c_{\ell +1}-c_{\ell })/2 + (c_{\ell }-c_{\ell -1})/2,\quad\text{for }\ell =2,\ldots M_{[a,b]}-1,\notag\\
	&w_1 = (c_{2}-c_{1})/2 + (c_{1}-a),\notag\\
	&w_{M_{[a,b]}} = (b-c_{M_{[a,b]}})+(c_{M_{[a,b]}}-c_{M_{[a,b]}-1})/2.
	\label{eq:weights}
\end{align}
We motivate this choice of weights in \S\ref{sec:SVD}, just after \rf{richness}.

\section{Discussion about the SVD solver}
\label{sec:SVD}

In this section we elaborate on the motivation for considering oversampling in combination with a truncated SVD in the solution method \eqref{eq:solution_pinv}. It is known from the earlier analysis in \cite{hewett2014frequency} -- recall Theorem \ref{dudnThm} in this paper -- that the best approximation in the HNA space $\VHNA$ converges to $v-\Psi$ 
at an exponential rate. Briefly, oversampling and regularization ensure that a near-best approximation can also be computed numerically, in spite of potential ill-conditioning of the linear system~\rf{eq:systemDeets}.

One obvious source of ill-conditioning of the discretization is that the basis functions $\varphi_N$ of the approximation space \eqref{HNASpace} may be close to being linearly dependent. The underlying reason is that, loosely speaking, in our discretization on each segment $\Gamma_j$ we are combining two bases together. The impact of this potential redundancy depends on the values of the wavenumber $k$ and the degree vector $\mathbf{p}$. Consider, for example, the case of small $k$: then the functions in the space $\calP_{\mathbf{p}}(\calM_j^+)\re^{\ri k s}$ are smooth and non-oscillatory on $\Gamma_j$, and they may approximate the functions in the space $\calP_{\mathbf{p}}(\calM_j^-)\re^{-\ri k s}$. In the case of fixed $k$ and a degree vector $\mathbf{p}$ associated with a large maximum degree $\pmax$ in \eqref{piDef1}, we can make a similar observation: the large degree polynomials in $\calP_{\mathbf{p}}(\calM_j^\pm)$ may resolve the oscillations of $\re^{\pm \ri k s}$. 

For linear systems with a numerical null-space, the truncated SVD solution satisfies a specific algebraic minimization property:
\begin{lemma}[\cite{coppe2019AZ}, Lemma 3.1]\label{lem:svd}
	Let $\bfa$ be computed by the regularized pseudo-inverse \eqref{eq:solution_pinv} with relative threshold $\epsilon > 0$. Then
	\begin{align}
		\label{eq:svd}
		\| \bfb - A \bfa \|_2 \leq \inf_{\bfv \in \mathbb{C}^N} \left\{ \| \bfb - A \bfv\|_2 + \epsilon \,\|A\|_2\,\| \bfv \|_2 \right\}.
	\end{align}
\end{lemma}
Lemma~\ref{lem:svd} shows that the regularized pseudo-inverse \eqref{eq:solution_pinv} yields a solution with small residual, provided that a coefficient vector $\bfv$ exists with small residual \emph{and} sufficiently small norm $\| \bfv \|_2$. Note that the size of the coefficients is affected by the normalization of the basis functions, and recall from \S\ref{sec:HNAspace} that for our problem we have normalized the basis functions in $L^2(\Gamma)$.
For us, the existence of suitable vectors that result in a small right-hand side in \rf{eq:svd} is suggested by Theorem~\ref{vpmThm} and Theorem~\ref{dudnThm}. Thus, ill-conditioning of $A$ due to redundancy in the discretization does not preclude the computation of highly accurate solutions.

However, the statement of Lemma~\ref{lem:svd} is purely algebraic, and a small discrete Euclidean residual of \rf{eq:systemDeets} does not imply a small continuous residual of the integral equation \eqref{eq:HNABIE} in $H^{1/2}(\Gamma)$, from which we could infer (by the bounded invertibility of $S:\tilde H^{-1/2}(\Gamma) \to H^{1/2}(\Gamma)$) a small approximation error of the continuous solution of \eqref{eq:HNABIE} in $\tilde H^{-1/2}(\Gamma)$. 

In the simpler $L^2$ setting one can appeal to the results on function approximations using redundant sets and a regularizing SVD solver with threshold $\epsilon$ presented in \cite{AdHu:19:FNA,AdHu:18:FNA2}. No guarantees can be given about the accuracy of interpolation, with errors possibly as large as $1/\epsilon$ \cite[Proposition 4.6]{AdHu:18:FNA2}. However, oversampling by a `sufficient' amount renders the approximation problem well-conditioned \cite[\S 6]{AdHu:19:FNA}. 
In particular, if functions in a space $G$ are sampled using a family of functionals $\{ \ell_{M,m} \}_{m=1}^M$, with $\ell_{M,m} : G \to \mathbb{C}$, then ideally the samples should be sufficiently `rich' to recover any function in $G$:
\begin{equation}\label{richness}
	A' \| g \|^2 \leq \liminf_{M \rightarrow \infty} \sum^{M}_{m=1} | \ell_{m,M}(g) |^2 \leq \limsup_{M \rightarrow \infty} \sum^{M}_{m=1} | \ell_{m,M}(g) |^2 \leq B' \| g \|^2,\quad \forall g \in G.
\end{equation}
The choice to weight matrix entries with $\sqrt{w}_m$ in \eqref{eq:systemDeets}, and the specific choice \eqref{eq:weights} of the weights, yields $A'=B'=1$ for $G \subset L^2(\Gamma)$, because the discrete sums are a Riemann sum for the $L^2$ norm.

Generalising these results to the $\tilde{H}^{-1/2}-H^{1/2}$ setting required for the rigorous analysis of our collocation BEM remains an open problem. Therefore, aside from the motivating observations presented in this section, our choice of collocation points, and of a linear oversampling rate $M \approx C_{\text{OS}} N$ with a proportionality constant $C_{\text{OS}}$ that is independent of the wavenumber, is based largely on numerical experiments, results of which we report in \S\ref{sec:numResults} (and see also \cite{emile}, where a number of different collocation point allocations were compared for square systems without oversampling).

\section{Oscillatory quadrature}\label{sec:oscillatoryquadrature}

Our HNA approximation space uses oscillatory basis functions supported on large ($k$-independent) mesh elements. This means that assembling the discrete system \eqref{eq:systemDeets} involves calculating highly oscillatory integrals, which may also be singular.

By applying linear changes of variable and possibly splitting the integration range, 
the integrals arising in \eqref{eq:systemDeets} can all be written in terms of integrals of the general form
\begin{align}
	\label{eq:GeneralIntegral}
	I[F;k,T,\alpha] = \int_0^T F(t; k)\re^{\imag k \alpha t}\rd t,
\end{align}
where $k>0$ is the wavenumber, $T>0$, $\alpha\in[0,2]$, and $F(\cdot;k)$ is smooth and non-oscillatory on $(0,T)$ but possibly logarithmically singular at, or near, the left endpoint $t=0$. More explicitly, we shall assume that, for some $t_0\in(-\infty,0]$, $F(\cdot;k)$ is analytic in the cut-plane $\C\setminus(-\infty,t_0]$, with at most polynomial growth as $|t|\rightarrow\infty$, and that there exists $\tilde{c}>0$ and two functions $F_0(t;k)$ and $F_1(t;k)$, analytic in $k|t-t_0|< \tilde{c}$, such that 
\begin{equation}\label{GGQform}
	F(t;k)= F_0(t;k) + F_1(t;k)\log(k|t-t_0|), \quad\text{for } k|t-t_0|< \tilde{c}.
\end{equation} 
We explain the transformation of the integrals in \eqref{eq:systemDeets} to the form \eqref{eq:GeneralIntegral} in \S\ref{sec:transformation}. Then in \S\ref{sec:NSD} we outline our method for efficiently calculating integrals of the form \eqref{eq:GeneralIntegral}. 

\subsection{Transformation to the general form \eqref{eq:GeneralIntegral}}\label{sec:transformation}

The left-hand side of the system \eqref{eq:systemDeets} consists of integrals of the form
\begin{align}
	\label{eq:matrixNSD}
	S\varphi_n(c_m)=&{\frac{\imag}{4}\int_{a}^{b}H^{(1)}_0(k|c_m-s|)}{\calL_{q,[a,b]}(s)\re^{\pm\imag k s }}\rd{s},
	\nonumber\\
	=&\frac{\imag}{4}\int_{a}^{b}{\frac{ H^{(1)}_0(k|c_m-s|)}{\re^{\imag k|c_m-s|}}}{\calL_{q,[a,b]}(s)}\re^{\imag k (|c_m-s|\pm s)}\rd{s},
\end{align}
where $\varphi_n=\calL_{q,[a,b]}(s)\re^{\pm\imag k s }$ is the $n$th basis function, 
$\calL_{q,[a,b]}$ 
is the 
$L^2$-normalised 
Legendre polynomial %\eqref{eq:LegRS} 
of some degree $q$ on $[a,b]=\supp{\varphi_n}$,  
and $c_m$ is the $m$th collocation point. 
To obtain a non-oscillatory prefactor we have extracted the phase of the Hankel function in \rf{eq:matrixNSD} using \cite[(10.2.5)]{NIST:DLMF}. 
To express \rf{eq:matrixNSD} in the form \rf{eq:GeneralIntegral} we then proceed by cases. In each case below, the fact that the resulting prefactor satisfies \eqref{GGQform} follows from the small argument behaviour of the Hankel function (see e.g. \cite[(10.4.3),(10.2.2),(10.8.2)]{NIST:DLMF}). 

\begin{itemize}
	\item
	\underline{Case A}: 
	If $c_m\leq a$ then $|c_m-s|=s-c_m$ on $[a,b]$ and a translation $s=a+t$ puts \rf{eq:matrixNSD} in the form \rf{eq:GeneralIntegral} with $T=b-a$, $\alpha=1\pm 1$ (i.e.\ $0$ for $+$ and $2$ for $-$) and \[F(t;k) = \frac{\imag}{4}\re^{\imag k(a\pm a-c_m)}\frac{ H^{(1)}_0(k(a+t-c_m))}{\re^{\imag k(a+t-c_m)}}{\calL_{q,[a,b]}(a+t)},\] which has a logarithmic singularity at $t_0 = c_m-a\leq 0$.
	
	\item
	\underline{Case B}: 
	If $c_m\geq b$ then the reflection $t\mapsto -t$ puts us in Case A.
	
	\item
	\underline{Case C}: 
	If $a<c_m<b$ then splitting the integral as $\int_a^b = \int_a^{c_m} + \int_{c_m}^b$ produces two integrals satisfying the conditions of Cases B and A respectively.
\end{itemize}

The right-hand side of \eqref{eq:LS} involves the evaluation of 
\begin{equation}\label{eq:oscStyle}
	S \Psi(c_m)=\frac{\resub{|d_2|}k}{2}\sum_{j=1}^{\numScreens}\int_{s_{2j-1}}^{s_{2j}}\frac{ H^{(1)}_0(k|c_m-s|)}{\re^{\imag k|c_m-s|}}\re^{\imag k (|c_m-s|+d_1s)}\rd{s}.
\end{equation}
Each integral in the sum \rf{eq:oscStyle} can be expressed in the form \rf{eq:GeneralIntegral} by essentially the same procedure described above. For example, when $c_m\leq s_{2j-1}$ we can adapt the approach of Case A to write $\int_{s_{2j-1}}^{s_{2j}}$ in the form \rf{eq:GeneralIntegral} with 
$T=L_j=s_{2j}-s_{2j-1}$, $\alpha=1+d_1$, and $F(t;k) = \frac{\resub{|d_2|}k}{2}\re^{\imag k(s_{2j-1}-c_m)}\frac{ H^{(1)}_0(k(s_{2j-1}+t-c_m))}{\re^{\imag k(s_{2j-1}+t-c_m)}}$, which has a logarithmic singularity at $t_0 = c_m-s_{2j-1}\leq 0$.

\subsection{Fast evaluation of \rf{eq:GeneralIntegral} using numerical steepest descent with generalised Gaussian quadrature}
\label{sec:NSD}

Depending on the values of the parameters $k,T,\alpha$ and $t_0$, the integrand \rf{eq:GeneralIntegral} may be highly oscillatory and/or numerically singular. In this section we describe an algorithm which can compute \rf{eq:GeneralIntegral} efficiently in all cases. 
Our approach 
combines the method of numerical steepest descent (see e.g.\ \cite{DeHu:09,HuVa:06,DaanBook}), which uses complex contour deformation to convert oscillatory integrals into rapidly decaying ones, with generalised Gaussian quadrature (see e.g.\ \cite{HuCo:08,huybrechs2017computation}), which handles the singularities.
Our method is an extension of the scheme presented in \cite[\S4.3]{huybrechs2017computation}, which was shown to accurately calculate singular oscillatory integrals, to the case of near-singularities, where the integrand has a singularity outside of, but close to the integration range. 

Our algorithm involves two parameters $c_{\rm osc}>0$ and $0<c_{\rm sing}\leq 1$, which are used to classify the integral \rf{eq:GeneralIntegral} into one of four cases 
(described below). The parameter $c_{\rm osc}$ represents the minimum number of oscillations in the exponential factor $\re^{\ri k \alpha t}$ there need to be over the interval $[0,T]$ for us to consider the integral oscillatory (and apply numerical steepest descent). That is, we consider the integral oscillatory when $k\alpha T/(2\pi)>c_{\rm osc}$ and non-oscillatory otherwise. 
The parameter $c_{\rm sing}$ controls how close to the integration range the singularity at $t_0$ needs to be for us to consider the integral singular (and apply generalised Gaussian quadrature). An important factor affecting the accuracy of numerical quadrature based on polynomial approximation for non-entire functions is the distance from the integration interval to the nearest singularity, measured \emph{relative to the length of the integration interval} (e.g. \cite[Theorem~19.3]{Trefethen2013}). So, in the non-oscillatory case we consider the integral singular when $|t_0|/T<c_{\rm sing}$. But when the integral is oscillatory this is an unnecessarily stringent condition, because as $k\alpha\to \infty$ with $T$ fixed the main contribution to the integral comes from small neighbourhoods of the endpoints $t=0$ and $t=T$ of size approximately $O(1/(k\alpha))$ (e.g. \cite[(2.1)]{DaanBook}). Hence, in the oscillatory case we consider the integral singular when $|t_0|k\alpha /(2\pi c_{\rm osc}) <c_{\rm sing}$. The inclusion of $c_{\rm osc}$ in the denominator here ensures our  classifications are compatible, in the sense that in the borderline case $k\alpha T/(2\pi)=c_{\rm osc}$, the transition between singular and non-singular cases occurs at $|t_0|/T = {c}_{\rm sing}$ for both the oscillatory and non-oscillatory cases.
%$|t_0|k\alpha /(2\pi) <\tilde{c}_{\rm sing}$. 

Our algorithm also depends on one further parameter, $\numQuad\in\N$, which controls the number of quadrature points used. 
Explicitly, we use $\numQuad$ points for each standard Gaussian quadrature rule and $2\numQuad$ points for each generalised Gaussian quadrature rule (since generalised Gaussian quadrature requires twice as many points as standard Gaussian quadrature to achieve the same degree of exactness, see e.g.\ \cite{HuCo:08}). 

We now present our algorithm, along with diagrams illustrating the contour deformations involved. Here the arrows indicate the direction along which each contour should be traversed, and the dashed line represents the original integration contour $[0,T]$.

\begin{itemize}
	\item 
	
	\noindent\underline{\textbf{Case 1}}: $k\alpha T/(2\pi) \leq c_{\rm osc}$ and $|t_0|/T\geq c_{\rm sing}$ (non-oscillatory, non-singular). 
	We evaluate \rf{eq:GeneralIntegral} using standard Gauss-Legendre quadrature with $\numQuad$ points. 
	
	\item 
	\noindent\underline{\textbf{Case 2}}: $k\alpha T/(2\pi) \leq c_{\rm osc}$ and $|t_0|/T< c_{\rm sing}$ (non-oscillatory, singular). 
	We write \rf{eq:GeneralIntegral} as the difference of two singular integrals
	
	\begin{tabular}{b{0.4\linewidth} b{0.5\linewidth}}
		\[ \int_{0}^T =- \int_{t_0}^0+ \int_{t_0}^T , \] 
		&
		%\vspace*{.0cm}
		\begin{tikzpicture}
		\coordinate (O) at (0,0);
		\coordinate (Om) at (0,-0.1);
		\coordinate (Op) at (0,0.1);
		\coordinate (T) at (3,0);
		\coordinate (t0) at (-.5,0);
		\coordinate (Om) at (0,-0.1);
		\coordinate (Tm) at (3,-0.1);
		\coordinate (t0m) at (-.5,-0.1);
		\coordinate (Op) at (0,0.1);
		\coordinate (Tp) at (3,0.1);
		\coordinate (t0p) at (-.5,0.1);
		\draw[dashed] (O)--(T);
		\draw[->] (t0m)--(Tm) ;
		\draw[->] (Op)--(t0p) ;
		%\filldraw 
		\filldraw (O) circle (1pt) node[anchor=south] {$0$};
		\filldraw (t0) circle (1pt) node[anchor=east] {$t_0$};
		\filldraw (T) circle (1pt) node[anchor=west] {$T$};
		\end{tikzpicture}
	\end{tabular}
	
	\noindent to which we apply generalised Gauss-Legendre quadrature with $2\numQuad$ points for each integral. (When $t_0=0$ the first integral is not present.) 
	\item 
	\noindent\underline{\textbf{Case 3}}: $k\alpha T/(2\pi) > c_{\rm osc}$ and $|t_0|k\alpha /(2\pi c_{\rm osc}) \geq c_{\rm sing}$ (oscillatory, non-singular). 
	We deform the integration contour into the complex plane, writing
	
	\begin{tabular}{b{0.4\linewidth} b{0.5\linewidth}}
		%		\vspace*{1cm}
		\[ \int_{0}^T = \int_{0}^{\ri \infty} - \int_{T}^{T+\ri \infty}, \]
		&%\vspace{0cm}
		\begin{tikzpicture}
		\coordinate (O) at (0,0);
		\coordinate (T) at (3,0);
		\coordinate (t0) at (-2,0);
		\coordinate (Oinf) at (0,2);
		\coordinate (Tinf) at (3,2);
		
		\draw[->] (O)--(Oinf);
		\draw[->] (Tinf)--(T) ;
		\draw[dashed] (O)--(T) ;
		%\filldraw 
		\filldraw (O) circle (1pt) node[anchor=north] {$0$};
		\filldraw (t0) circle (1pt) node[anchor=north] {$t_0$};
		\filldraw (T) circle (1pt) node[anchor=north] {$T$};
		\filldraw (Tinf)  node[anchor=south] {$T+\imag\infty$};
		\filldraw (Oinf)  node[anchor=south] {$\imag\infty$};
		\end{tikzpicture}
	\end{tabular}
	
	\noindent justified by our assumption that $k\alpha\geq 0$ and that $F(\cdot;k)$ grows only polynomially at infinity, meaning that the integrand of \rf{eq:GeneralIntegral} decays exponentially as ${\rm Im}\,t\to\infty$. Parametrizing the integrals by $t=\ri k\alpha \xi$ (respectively $t=T+\ri k\alpha \xi$) gives 
	\begin{equation}\label{eq:SD}
		\int_0^T = \frac{\ri}{k \alpha}\int_0^\infty F\left(\frac{\imag \xi}{k\alpha}; k\right)\re^{- \xi}\,\rd \xi - \frac{\ri\re^{\imag k\alpha T}}{{k\alpha}}\int_0^\infty F\left(T + \frac{\imag \xi}{\alpha k}; k\right)\re^{-\xi}\,\rd \xi;
	\end{equation}
	we evaluate both semi-infinite integrals using standard Gauss-Laguerre quadrature with $\numQuad$ points for each integral.
	
	\item 
	\noindent\underline{\textbf{Case 4}}: $k\alpha T/(2\pi) > c_{\rm osc}$ and $|t_0|k\alpha /(2\pi c_{\rm osc})< c_{\rm sing}$ (oscillatory, singular). 
	We write the integral as 
	
	\begin{tabular}{b{0.4\linewidth} b{0.5\linewidth}}
		%\vspace*{-2.7cm}
		\[ \int_{0}^T = -\int_{t_0}^0 + \int_{t_0}^{t_0+\ri \infty} - \int_{T}^{T+\ri \infty}, \]
		&
		\begin{tikzpicture}
		\coordinate (O) at (0,0);
		\coordinate (T) at (3,0);
		\coordinate (t0) at (-.8,0);
		
		\coordinate (t0inf) at (-.8,2);
		\coordinate (Oinf) at (0,2);
		\coordinate (Tinf) at (3,2);
		
		\draw[->] (O)--(t0);
		\draw[->] (t0)--(t0inf);
		\draw[->] (Tinf)--(T) ;
		\draw[dashed] (O)--(T) ;
		%\filldraw 
		\filldraw (O) circle (1pt) node[anchor=north] {$0$};
		\filldraw (t0) circle (1pt) node[anchor=north] {$t_0$};
		\filldraw (T) circle (1pt) node[anchor=north] {$T$};
		\filldraw (t0inf)  node[anchor=south] {$t_0+\imag\infty$};
		\filldraw (Tinf)  node[anchor=south] {$T+\imag\infty$};
		\end{tikzpicture}
	\end{tabular}
	
	\noindent evaluating the first integral using generalised Gauss-Legendre quadrature with $2\numQuad$ points, %(after parametrization), 
	the second using generalised Gauss-Laguerre quadrature\footnote{In related literature it is common for \emph{generalised Gauss-Laguerre} to refer to a quadrature rule for integrals of the form $\int_0^\infty x^{-\nu}f(x)\re^{-x}\rd{x}$, for $\nu>-1$ where $f$ is well-approximated by a polynomial, but we do not use this definition here. In this paper, we use \emph{generalised Gauss-Laguerre} to mean a quadrature rule which can efficiently evaluate $\int_0^\infty[\log(x)f(x)+g(x)]\re^{-x}\rd{x}$, where $f$ and $g$ are well-approximated by polynomials.} with $2\numQuad$ points, and the third using standard Gauss-Laguerre quadrature with $\numQuad$ points. 
	(Again, when $t_0=0$ the first integral is not present.)
	
	Here our assumption that $c_{\rm sing}\leq 1$ ensures that our treatment of the third integral $\int_{T}^{T+\ri \infty}$ as non-singular (being evaluated by standard, rather than generalised Gauss-Laguerre quadrature) is consistent with our classification of the two integrals in Case 3 as non-singular, in the sense that if $k\alpha T/(2\pi) > c_{\rm osc}$ then
	\[k\alpha(T-t_0)/(2\pi c_{\rm osc}) \geq k\alpha T/(2\pi c_{\rm osc}) > 1 \geq c_{\rm sing}.\]
\end{itemize}

The total number of quadrature points (a proxy for computational cost) 
is $\numQuad$ in Case 1, $4\numQuad$ in Case 2, $2\numQuad$ in Case 3, and $5\numQuad$ in Case 4 (excluding special cases where $t_0=0$). Therefore, to minimise computational cost, we should take $c_\mathrm{osc}$ as large as possible, and $c_\mathrm{sing}$ as small as possible, so as to be in the cheapest Case 1 (non-oscillatory and non-singular) as often as possible. But obviously one cannot take $c_\mathrm{osc}$ too large, or $c_\mathrm{sing}$  too small, else oscillations (respectively, singularities) will not be adequately resolved. We investigate the choice of these parameters in more detail in the next section.

\subsection{Choice of parameters}\label{sec:quadParamChoice}

To validate the quadrature scheme described in \S\ref{sec:NSD} and tune the parameters $c_\mathrm{osc}$, $c_\mathrm{sing}$ and $\numQuad$ we 
carried out a detailed set of experiments for the 
model integral
\begin{align}
	\label{eq:TestIntegral}
	\int_0^1 H_0^{(1)}(k(t-t_0)) \,\rd t,
\end{align}
which is of the form \eqref{eq:GeneralIntegral} with $\alpha=T=1$ and $F(t;k)=H_0^{(1)}(k(t-t_0))\re^{-\imag k t}$.
For all possible combinations of $c_\mathrm{osc}\in\{1,\ldots,5\}$, $c_\mathrm{sing}\in\{0,0.1,0.2,\ldots,1\}$ and $\numQuad\in\{5,10,15,20\}$, we computed the integral \rf{eq:TestIntegral} using our algorithm at a large number of values of the parameters $k$ and $t_0$ in the range $k\in(0,100]$, $t_0\in[-1,0]$. For each set of parameters we measured the relative error compared to a reference value for \rf{eq:TestIntegral}, computed using a composite Gauss rule with a large number of quadrature points and mesh grading to handle (near) singularities. Based on numerical experiments we believe this reference solution to be accurate to at least $14$ digits. 

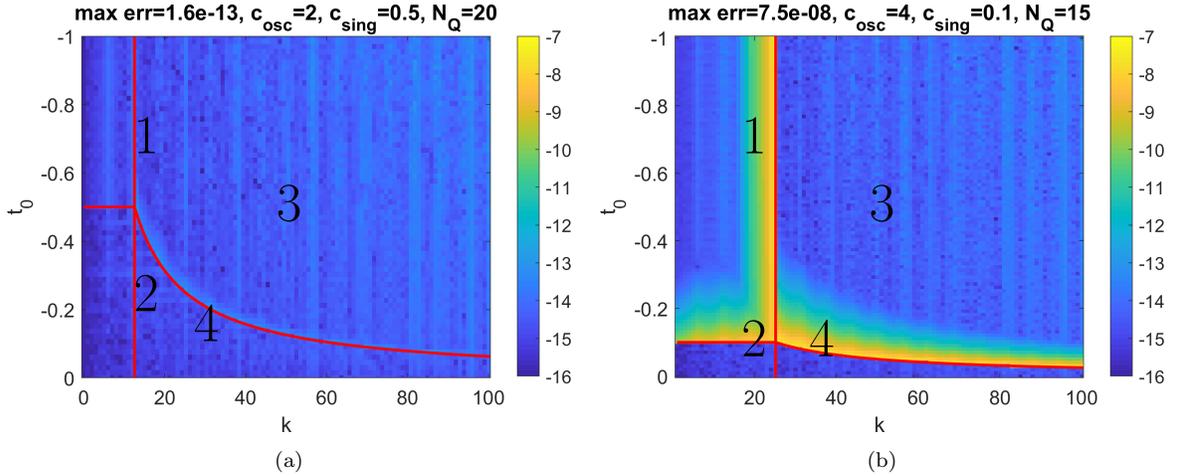
\begin{figure}%[ht]
	%\hspace{-1mm}
	\subfloat[]{
		\begin{tikzpicture}	
		\node at (0,0) {\includegraphics[width=.48\linewidth]{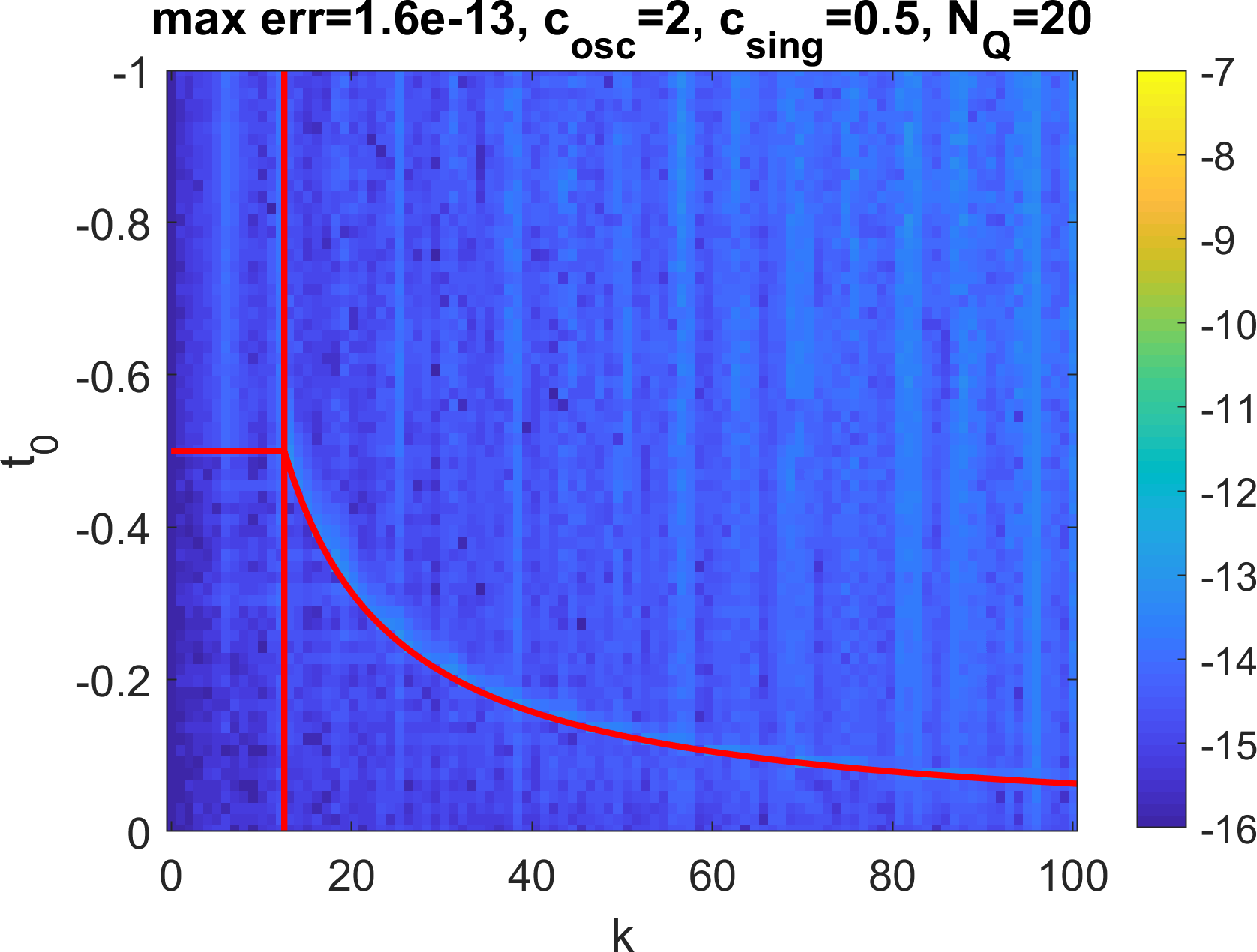}};
		\node at (-0,0.2) {\huge 3};
		\node at (-1.9,1.1) {\huge 1};
		\node at (-1.1,-1.4) {\huge 4};
		\node at (-1.9,-1) {\huge 2};
		\end{tikzpicture}
	}
	\subfloat[]{	
		\begin{tikzpicture}	
		\node at (0,0) {\includegraphics[width=.48\linewidth]{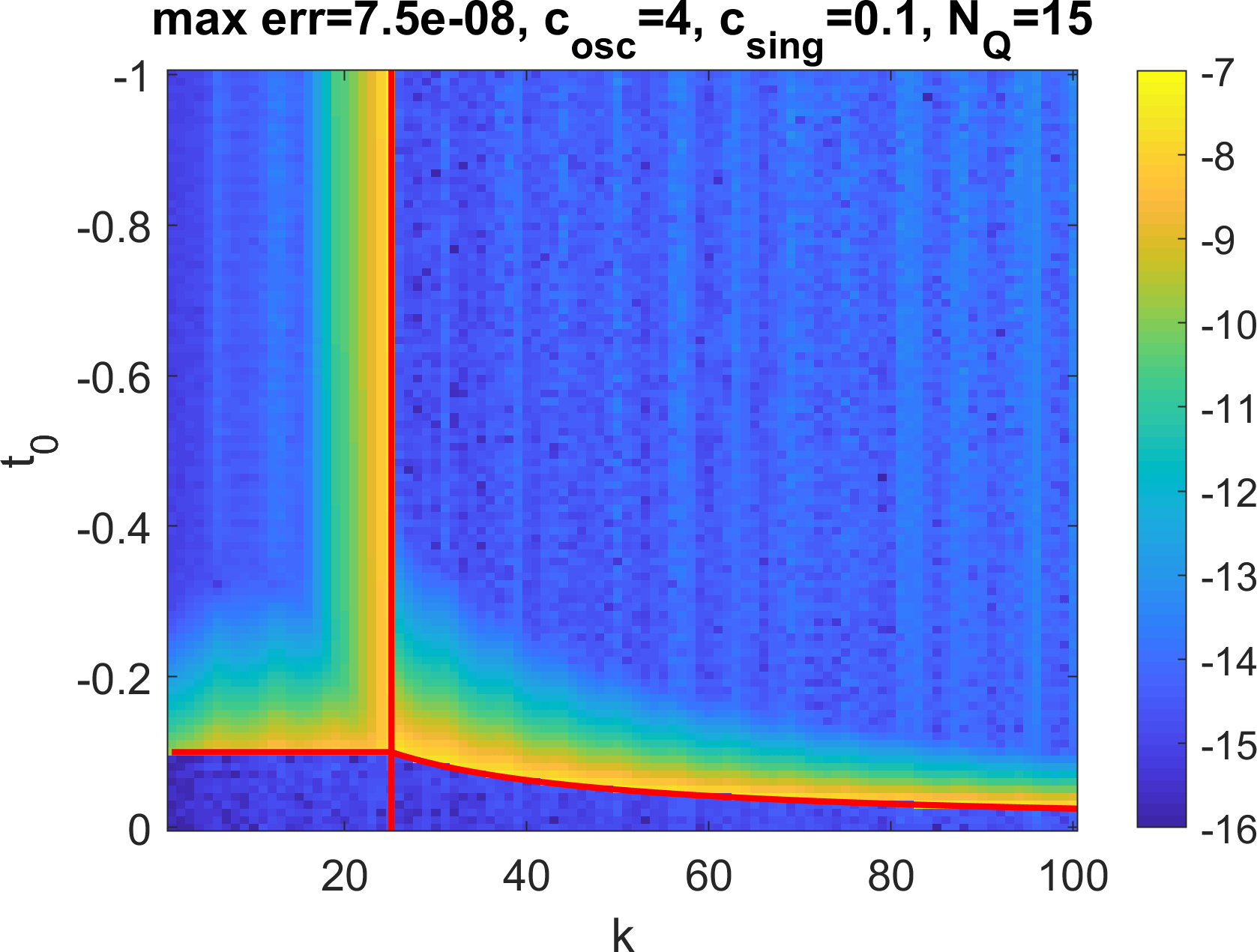}};
		\node at (-0,0.2) {\huge 3};
		\node at (-1.7,1.1) {\huge 1};
		\node at (-0.8,-1.6) {\huge 4};
		\node at (-1.7,-1.6) {\huge 2};
		\end{tikzpicture}	
	}	
	\caption{Plots of $\log_{10}$ of the relative error in computing \rf{eq:TestIntegral} using the quadrature scheme of \S\ref{sec:NSD} for two different sets of quadrature parameters $\{c_\mathrm{osc}=2, c_\mathrm{sing}=0.5, \numQuad=20\}$ and $\{c_\mathrm{osc}=4, c_\mathrm{sing}=0.1, \numQuad=15\}$. In both plots the labels 1-4 indicate the regions of the $(k,t_0)$ plane in which Cases 1, 2, 3 and 4 apply.
	}
	\label{fig:quadChoice}
\end{figure}

Based on our experiments, the choices 
\begin{align}
	\label{eq:quadparams}
	c_\mathrm{osc}=2, \quad c_\mathrm{sing}=0.5, \quad \numQuad=20
\end{align}
were found to produce a scheme which agrees with the reference solution to at least 12 digits, uniformly over the range of $k$ and $t_0$ tested, and these are the values we use in all our numerical results in \S\ref{sec:numResults}. 
As evidence of this claim we show in Fig.~\ref{fig:quadChoice}(a) a plot of the corresponding relative error over the studied range of $k$ and $t_0$. The 
$(k,t_0)$ plane is divided into the four regions in which Cases 1,2,3 and 4 apply. The maximum relative error over the whole plot is $1.6\times 10^{-13}$.

For comparison we show in Fig.~\ref{fig:quadChoice}(b) the corresponding plot for the choices $c_\mathrm{osc}=4$, $c_\mathrm{sing}=0.1$ and $\numQuad=15$. 
The maximum relative error over the whole plot is now $7.5\times10^{-8}$, and one can clearly see the accuracy degrading as one moves to the right in region 1 (increasing oscillation) and downwards in regions 1 and 3 (approaching singularity). 

As we shall show in \S\ref{sec:numResults}, the quadrature scheme in \S\ref{sec:NSD} with the parameter choices \rf{eq:quadparams} is sufficiently efficient to achieve our goal of frequency independent computational cost for the HNA collocation BEM. However, we do not claim that the quadrature scheme is completely optimal. In particular, we note that the error in numerical steepest descent (without singularities) behaves like $O((k\alpha)^{-2\numQuad-1})$ as $k\alpha\to\infty$ (see e.g.\ \cite[p90]{DaanBook}), so in Cases 3 and 4 savings could be made by reducing $\numQuad$ as $k\alpha$ grows. But we reserve such further optimisation for future work.

\section{Numerical results}
\label{sec:numResults}

In this section we demonstrate the computational efficiency of our HNA collocation BEM via a series of numerical examples. We also include the results of experiments exploring the influence of the parameters $C_{\rm OS}$ and $\epsilon$ in the 
truncated SVD solver. Finally we present an application of our method to the computation of high frequency scattering by high order prefractals of the middle-third Cantor set. 
The code used to produce the results in this section forms part of a larger Matlab-based HNA BEM software repository, which is downloadable from \texttt{github.com/AndrewGibbs/HNABEMLAB}. 

Throughout this section we denote by $v_p:=\phi_p+\Psi$ our HNA BEM approximation to the solution $v$ of the original BIE \eqref{BIE_sl}, where $\phi_p$ denotes our approximation of the solution $\phi$ of the BIE \eqref{eq:HNABIE}, the subscript $p$ indicating the maximum polynomial degree used in our HNA approximation space $V_N$, and $\Psi$ denotes the geometrical optics contribution. 
We denote by $v^\pm_{j,p}$ the corresponding numerical approximations to the amplitudes $v^\pm_j$ of the oscillatory functions $\re^{\pm \imag k s}$ in the decomposition \rf{Decomp}. 
From the boundary solution we obtain an approximation 
\[ u^s_p:=-\cS v_p\]
to the scattered field $u^s$ in $D$ (cf.\ \rf{eqn:D_RepThm}),   
and an approximation 
\begin{equation}\label{def:FarFieldapprox}
	u^\infty_p(\theta) :=-\int_{\tilde{\Gamma}}\re^{-\imag ks\cos\theta} v_p(s)\rd s
\end{equation} 
to the far-field pattern
\begin{equation}\label{def:FarField}
	u^\infty(\theta) :=-\int_{\tilde{\Gamma}}\re^{-\imag ks\cos\theta} v(s)\rd s,
\end{equation} 
which describes the far-field behaviour of $u^s$ via 
\[
u^s(\bfx)\sim u^\infty(\theta)\frac{\re^{\imag(kr +\pi/4)}}{2\sqrt{2\pi k r}},\quad\text{for }\bfx = r(\cos\theta,\sin\theta),\quad\text{as }\|\bfx\|_2=:r\rightarrow\infty.
\]

Unless otherwise stated, in our numerical results we use the following parameters:
\\\emph{HNA space parameters (\S\ref{sec:HNAspace}):}
\\We use the mesh grading parameter $\sigma=0.15$ and the number of layers $\ell=2(p+1)$ for each graded mesh, as in related studies (e.g. \cite{hewett2014frequency,ChGiLaMo:19}).
\\\emph{Collocation parameters (\S\ref{sec:collocation} and \S\ref{sec:SVD}):}
\\
We use the oversampling parameter $C_\mathrm{OS}=1.25$ and the SVD truncation parameter $\epsilon=10^{-8}$. These values were chosen based on the experiments summarised in \S\ref{sec:tuningParams}.
\\\emph{Quadrature parameters (\S\ref{sec:oscillatoryquadrature}):}
\\We take $c_\mathrm{osc}=2$, $c_\mathrm{sing}=0.5$ and $\numQuad=20$, as discussed in \S\ref{sec:quadParamChoice} (see \rf{eq:quadparams}).

\subsection{High frequency performance}

\begin{figure}
	\subfloat[]{\includegraphics[width=\linewidth]{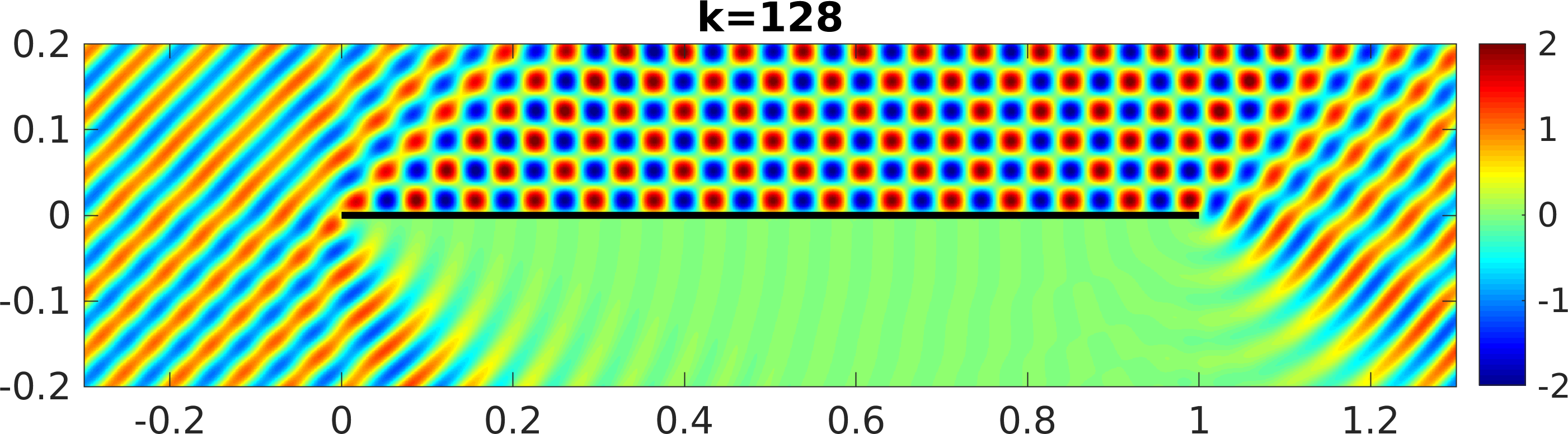}}\\
	\subfloat[]{\includegraphics[width=.48\linewidth]{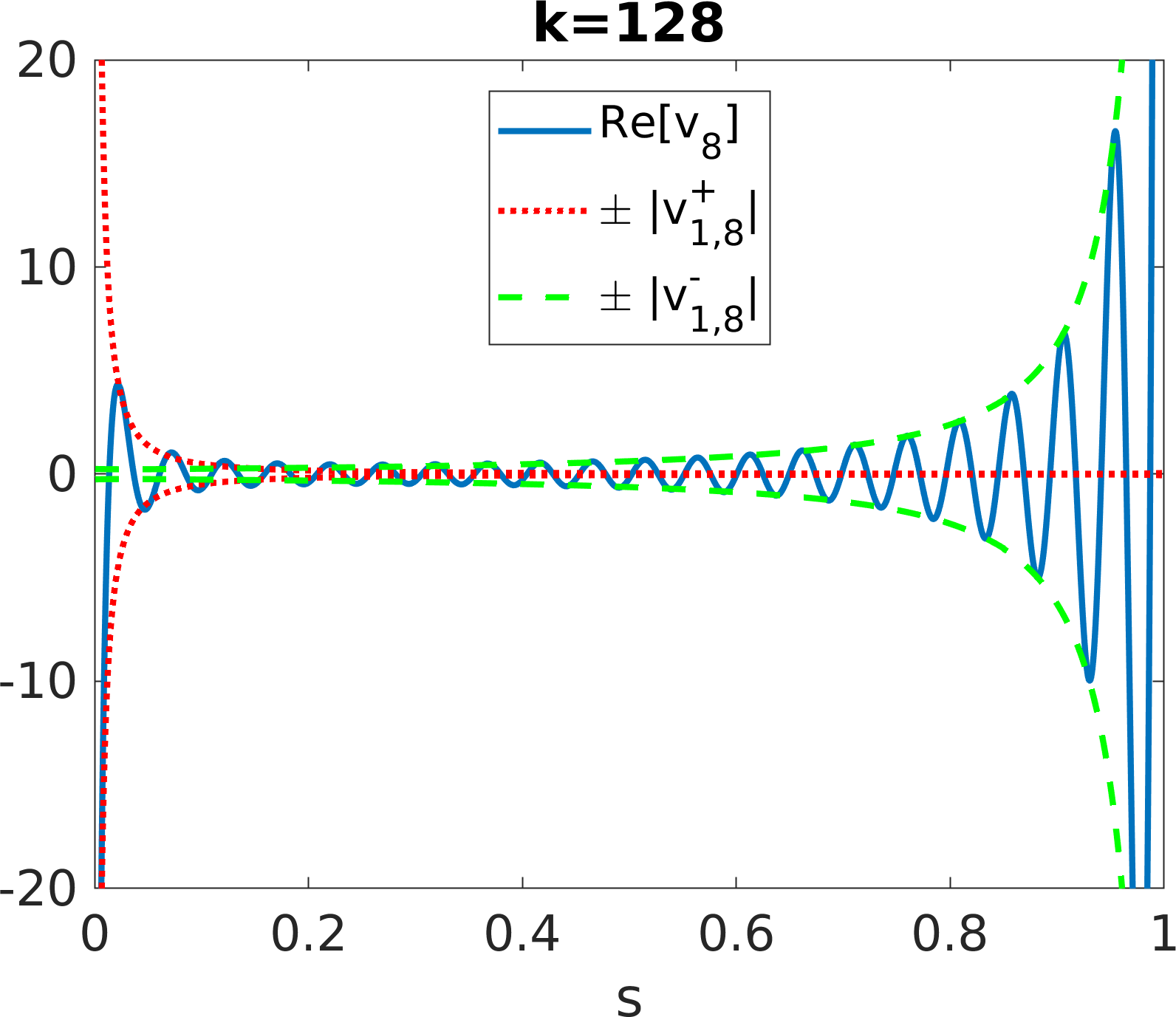}}
	\hspace{2mm}
	\subfloat[]{\includegraphics[width=.48\linewidth]{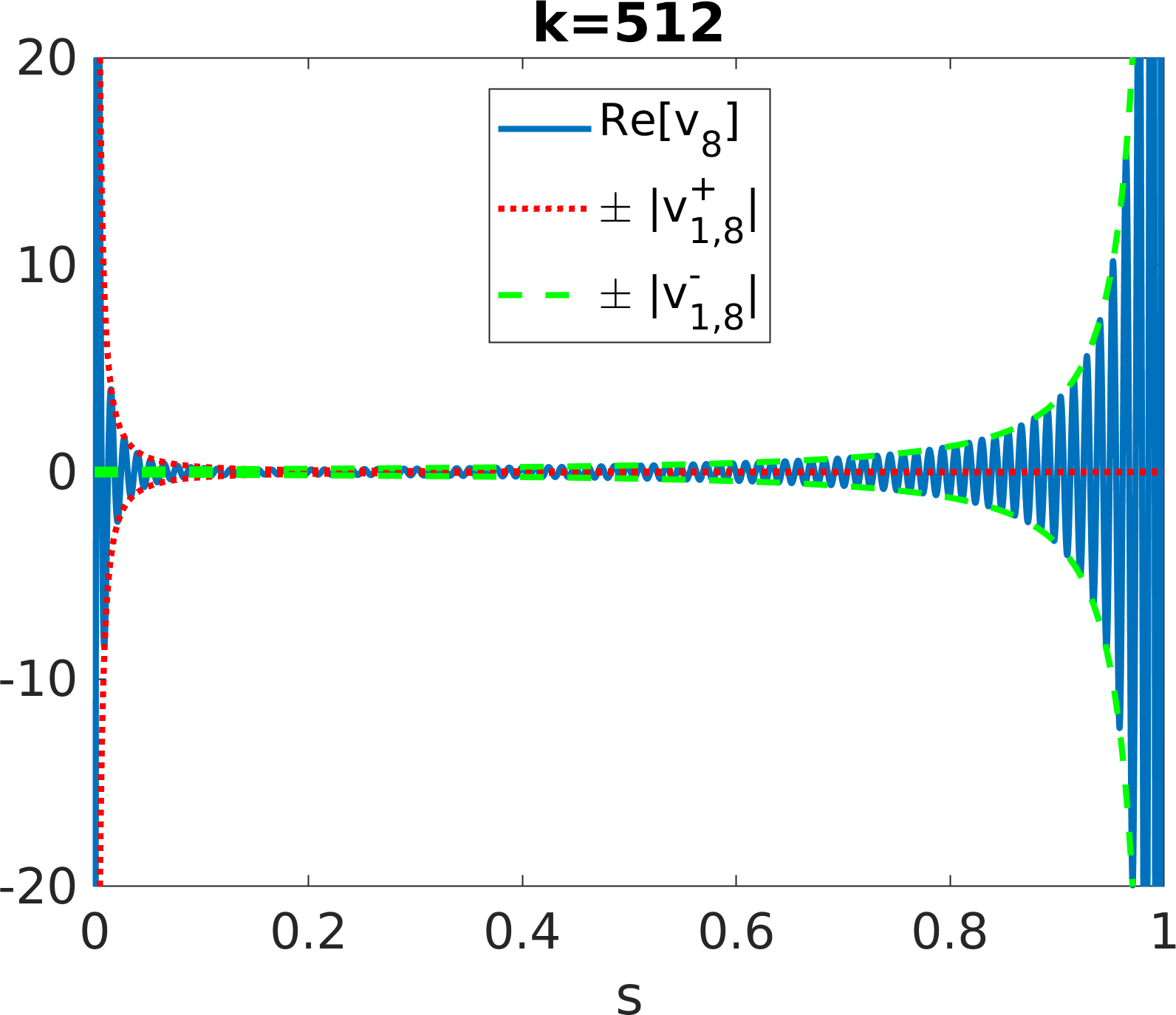}}
	\caption{
		HNA collocation BEM solution for $\Gamma=(0,1)\times\{0\}$ and $\bd =(1,-1)/\sqrt{2}$, with $p=8$. (a) Real part of the total field $u^i + u^s_8$ for $k=128$. (b) and (c): Real part of the boundary solution $v_8$, along with the amplitudes $|v^\pm_{1,8}|$, for $k=128$ and $512$ respectively.
	}
	\label{fig:referenceSoln}
\end{figure}

To illustrate the high frequency performance of our HNA method we consider first the case where the screen consists of a single unit interval $\Gamma=(0,1)\times\{0\}$. Examples with multiple screens will be considered in \S\ref{sec:Cantor}. As the incident wave direction we take $\bd =(1,-1)/\sqrt{2}$. Fig.~\ref{fig:referenceSoln}(a) shows the resulting total field for wavenumber $k=128$, and Figs.~\ref{fig:referenceSoln}(b) and \ref{fig:referenceSoln}(c) plot the boundary solution $v_8$, along with the magnitudes $|v^\pm_{1,8}|$ of the non-oscillatory amplitudes of the oscillatory factors $\re^{\pm \imag k s}$, for both $k=128$ and $k=512$. 

\begin{figure}%[ht]
	\subfloat[a][]{\includegraphics[width=.48\linewidth]{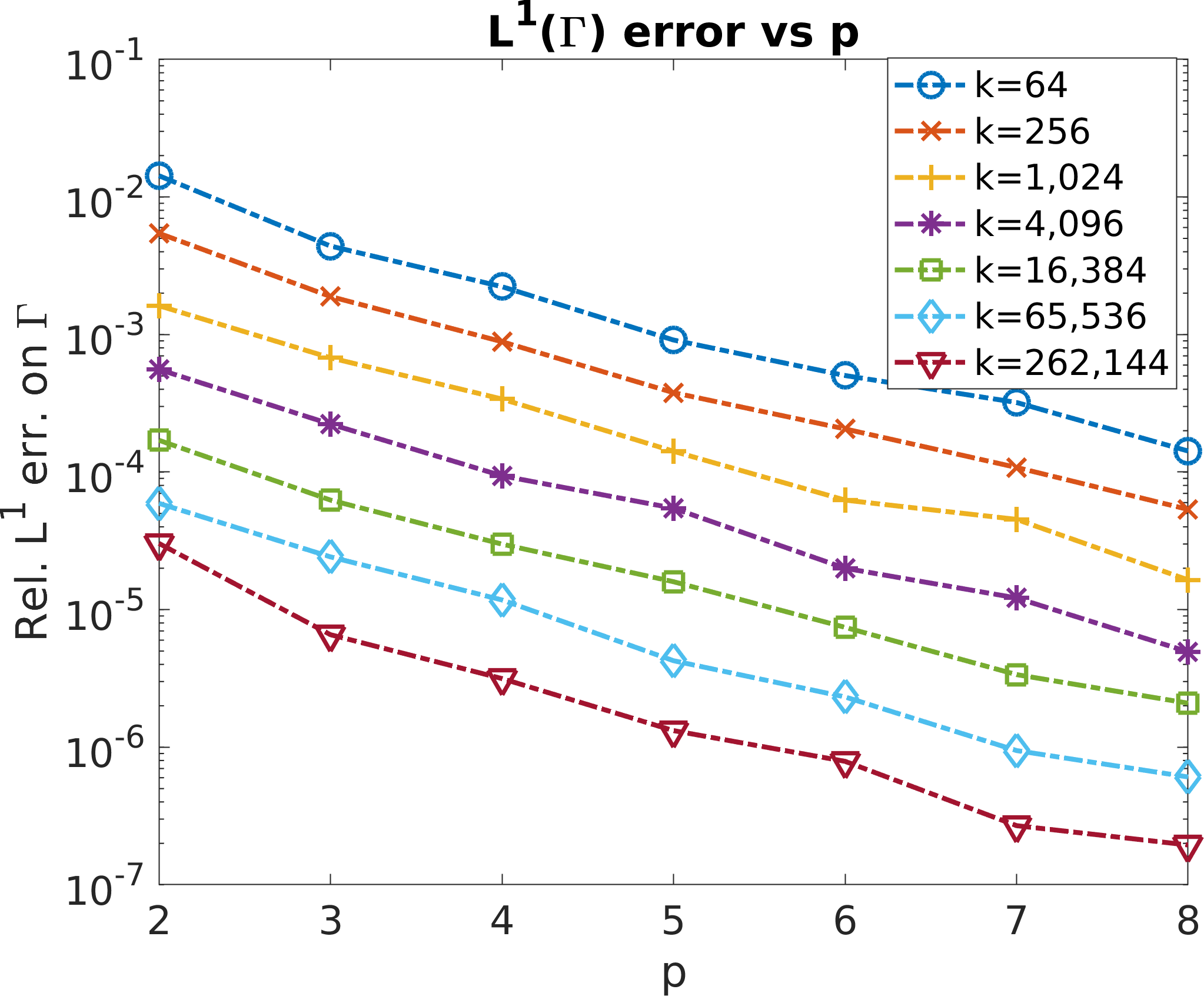}}\hspace{3mm}
	\subfloat[]{\includegraphics[width=.48\linewidth]{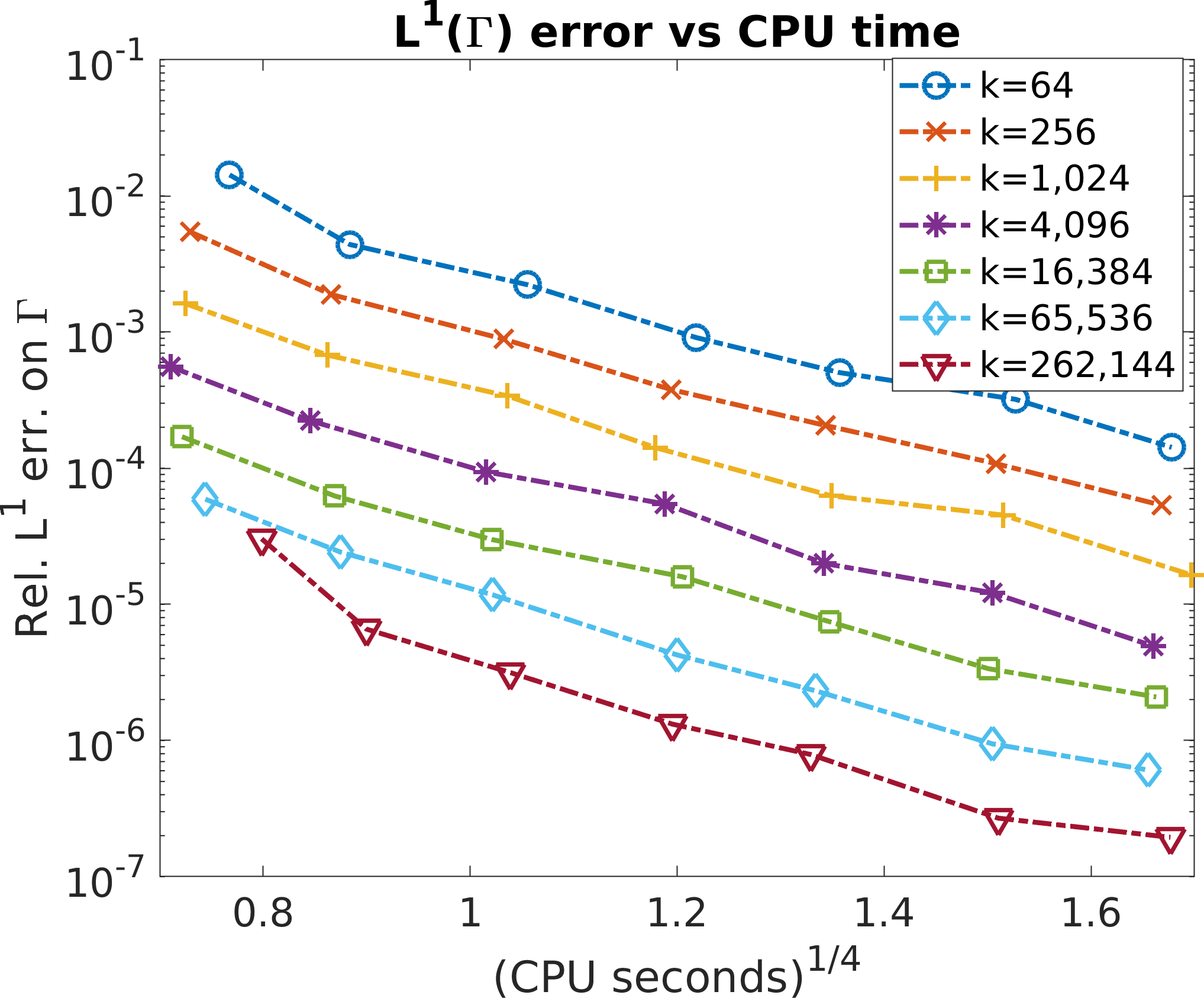}}\\
	\subfloat[]{\includegraphics[width=.48\linewidth]{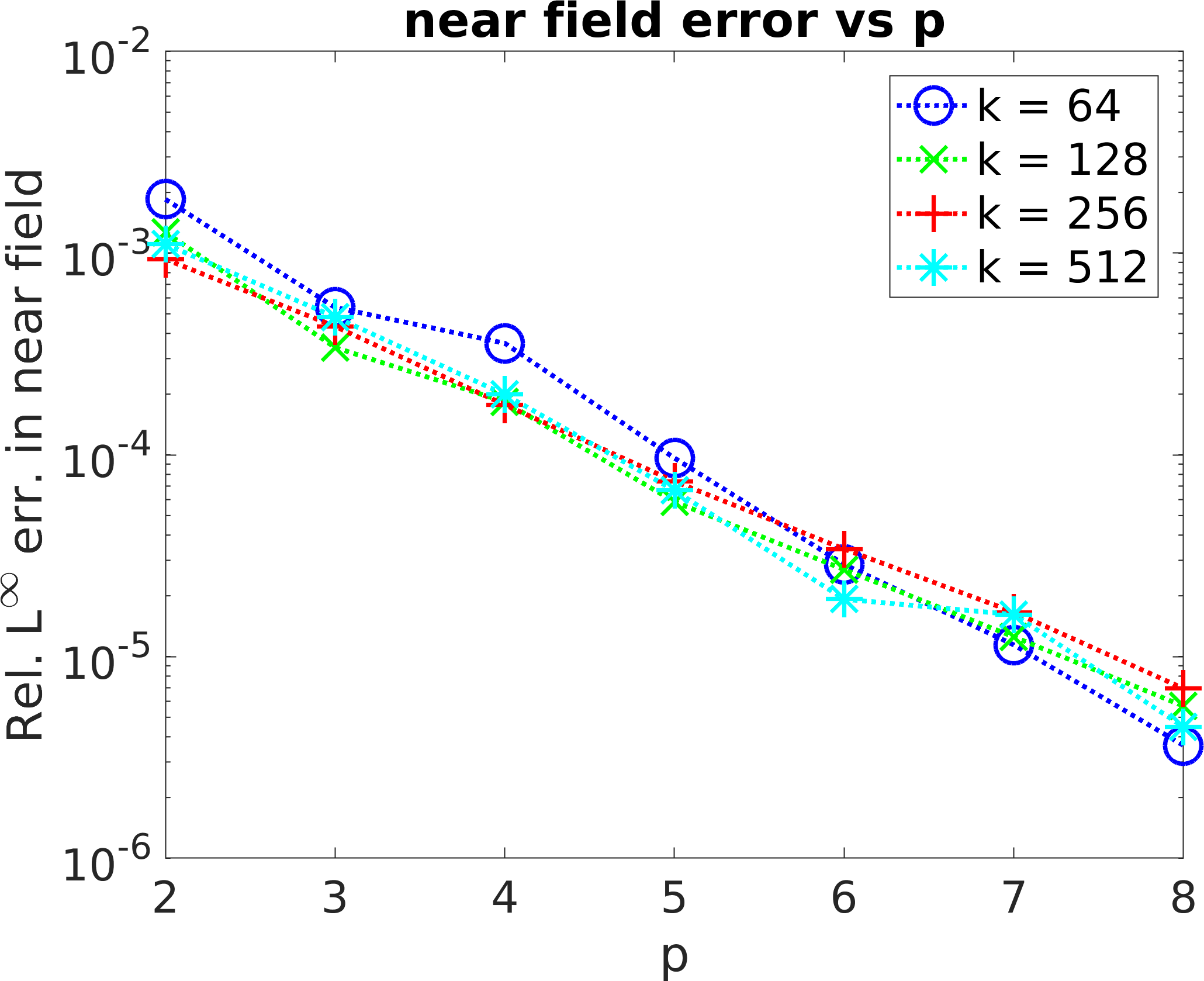}}\hspace{3mm}
	\subfloat[]{\includegraphics[width=.48\linewidth]{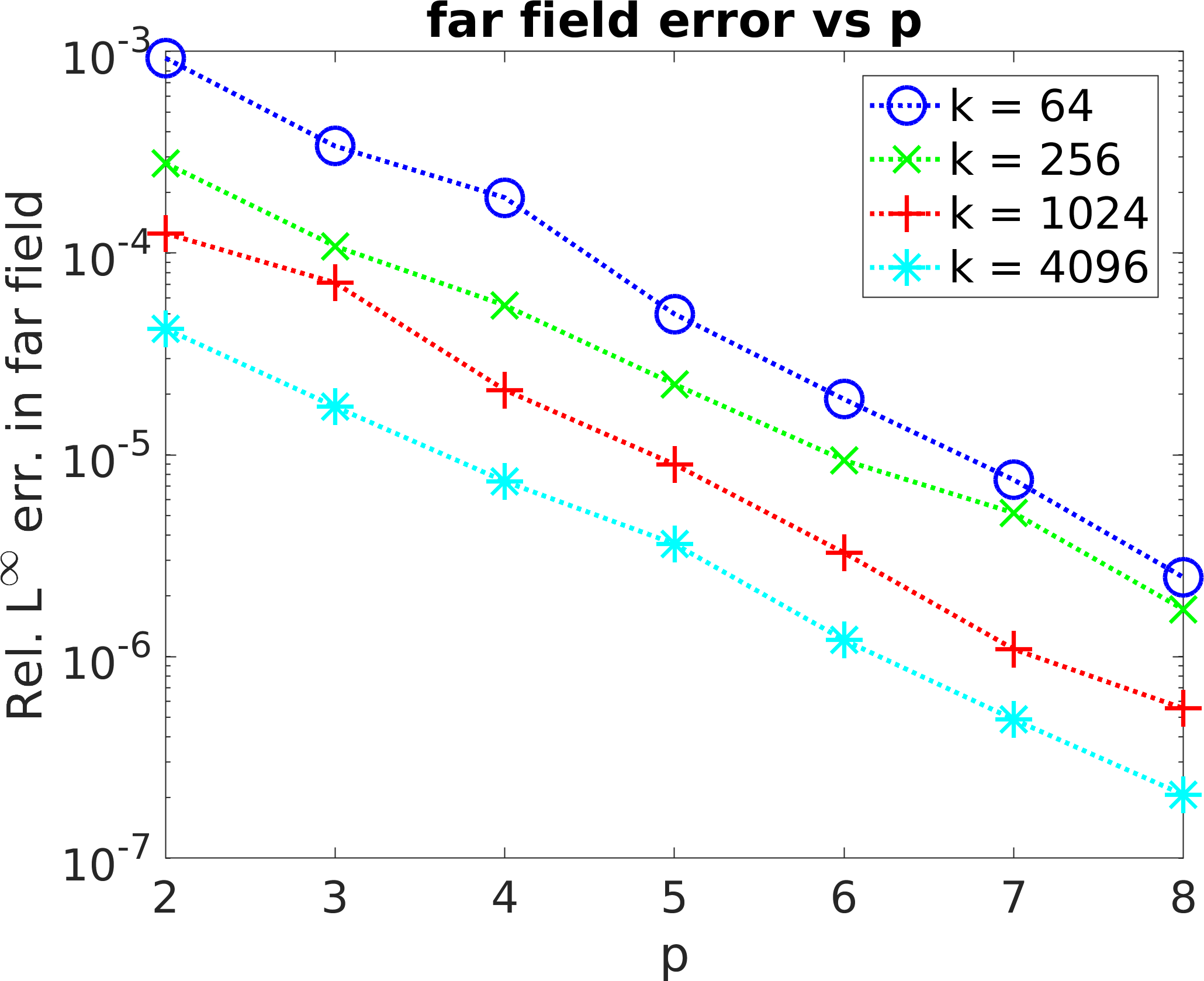}}
	\caption{Convergence with respect to increasing maximum polynomial degree $p$, for various wavenumbers $k$. (a) shows $L^1$ error on $\Gamma$ against $p$, (b) shows $L^1$ error on $\Gamma$ against CPU time, and (c) and (d) show $L^\infty$ error in the near and far field against $p$.
	}
	\label{fig:performance}
\end{figure}

In Fig.~\ref{fig:performance}(a) we show the relative $L^1(\Gamma)$ error 
\begin{equation}\label{eq:relErrL1Def}
	\text{Rel. }L^1\text{ err. on }\Gamma :=\frac{\|v_p-v_{12}\|_{L^1(\Gamma)}}{\|v_{12}\|_{L^1(\Gamma)}}
\end{equation}
in the BEM solution for a range of values of $p$ and $k$, using $p=12$ as reference solution. 
For easier comparison with the best approximation theory (Theorem \ref{dudnThm}) it would be preferable to compute the $\tilde H^{-1/2}(\Gamma)$ norm (which is possible via an appropriate single layer boundary integral representation, see e.g.\ \cite[\S7]{ChHeMoBe:19}), but here we choose the $L^1(\Gamma)$ norm because it is easier to evaluate, and also because convergence of the boundary solution in $L^1(\Gamma)$ implies, via a straightforward integral estimation, the convergence in supremum norm of the domain solution (on compact subsets of $D$) and the far-field pattern (cf.\ \rf{eqn:BdytoFF} below), which we study in Figs.~\ref{fig:performance}(c) and \ref{fig:performance}(d). %But we note that 
We compute our $L^1(\Gamma)$ norms by composite Gaussian quadrature with $20$ Gauss points per element on a graded mesh, built by intersecting the overlapping meshes used to approximate $v_{12}$, with additional subdivision of the larger elements to ensure that no element is longer than one wavelength $2\pi/k$.

Fig.~\ref{fig:performance}(a) demonstrates that our approximation of the boundary solution converges exponentially (in $L^1(\Gamma)$) as we increase the polynomial degree $p$. Moreover, for fixed $p$ the $L^1$ error is not only bounded but actually \textit{decreases} with increasing $k$. 

In Fig.~\ref{fig:performance}(b) we present the same results, plotted against the fourth root of the total CPU time for the BEM assembly and solve. The rationale behind this comparison is as follows. 
Our approximation space $V_N$ for this problem has dimension $N=O(p^2)$ (maximal polynomial degree $p$ over two overlapping meshes of $2(p+1)$ elements each). And with a fixed oversampling parameter $C_\mathrm{OS}=1.25$, the number of collocation points $M\approx C_\mathrm{OS}N$ also scales like $O(p^2)$. Hence the number of elements in the BEM matrix scales like $MN=O(p^4)$ with increasing $p$. The linear systems we use are moderate in size (typically $N\leq 1000$) and are rapidly solved using our truncated SVD approach. As a result, the majority of the CPU time is spent assembling the discrete system, which means that total CPU time scales like $O(p^4)$. Therefore, to scrutinize the performance of our numerical quadrature scheme for evaluating the BEM matrix entries by a direct comparison with Fig.~\ref{fig:performance}(a), we should plot the boundary error against the fourth root of the CPU time. The results in Fig.~\ref{fig:performance}(b) demonstrate that our quadrature scheme (based on numerical steepest descent) is both accurate and efficient, and results in a fully discrete BEM for which the error at a fixed computational cost is bounded (in fact, decreases) with increasing $k$. The timings in Fig.~\ref{fig:performance}(b) are for a desktop PC with an Intel 
i7-4790 3.60GHz 4-core CPU, with matrix assembly done in parallel using a Matlab \texttt{parfor} loop.
They show that highly accurate solutions can be obtained for essentially arbitrarily high frequencies in just a few seconds on a standard machine.

In Figs.~\ref{fig:performance}(c) and \ref{fig:performance}(d) we plot the corresponding near-field errors
\begin{equation}\label{eq:relErrNFDef}
	\text{Rel. }L^\infty\text{ err. in near-field} :=\frac{\|u^s_p-u^s_{12}\|_{L^\infty(\mathcal{O})}}{\|u^s_{12}\|_{L^\infty(\mathcal{O})}},
\end{equation}
where $\mathcal{O}$ is the circle of radius one centred at $(0.5,0)$, 
and far-field errors
\begin{equation}\label{eq:relErrFFDef}
	\text{Rel. }L^\infty\text{ err. in far-field} :=\frac{\|u^\infty_p-u^\infty_{12}\|_{L^\infty(0,2\pi)}}{\|u^\infty_{12}\|_{L^\infty(0,2\pi)}}.
\end{equation}

As for the boundary solution, we observe exponential convergence as $p$ increases with $k$ fixed, and bounded errors as $k$ increases for $p$ fixed. %\resub{In the far field the accuracy noticeably improves as $k$ increases. This is due to the fact that the far-field pattern $u^\infty$ has peaks of magnitude proportional to $d_2k$ at the values of $\theta$ corresponding to the reflected and shadow directions (cf.\ Fig.~\ref{fig:FFCantor}(a)), associated with the contribution of the geometrical optics term $\Psi$. Hence the denominator in the relative $L^\infty$ error will be $O(k)$ for non-grazing incidence, while the numerator is not expected to be comparably large because the contribution from $\Psi$ is computed exactly (up to quadrature error), and hence does not contribute to the discretization error.}

\begin{figure}[t]
	\subfloat[][]{\includegraphics[width=\linewidth]{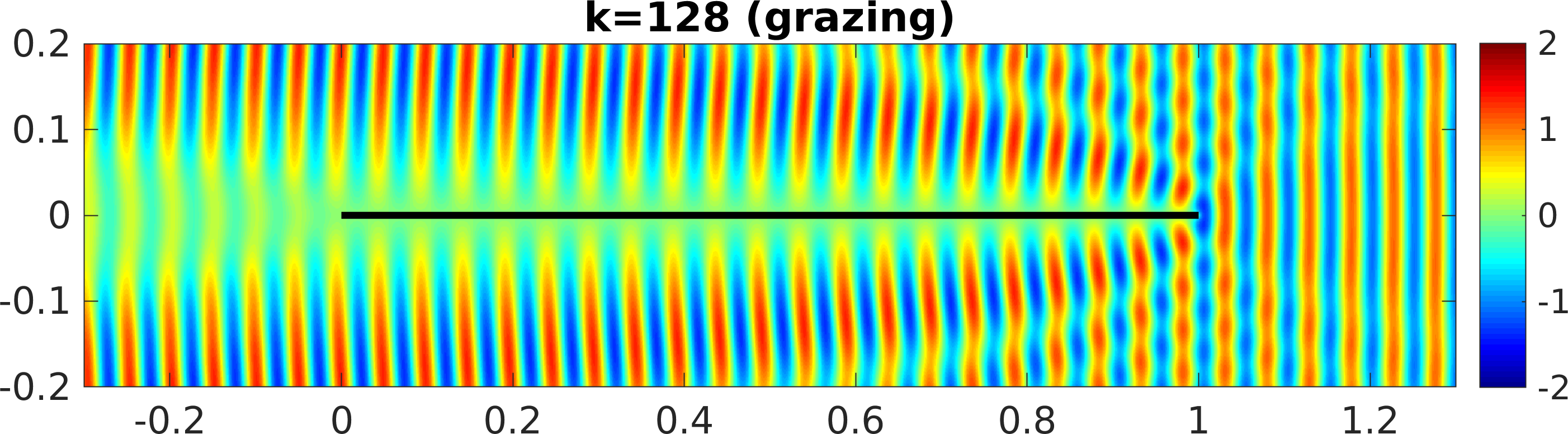}}\\
	\subfloat[][]{\includegraphics[width=.48\linewidth]{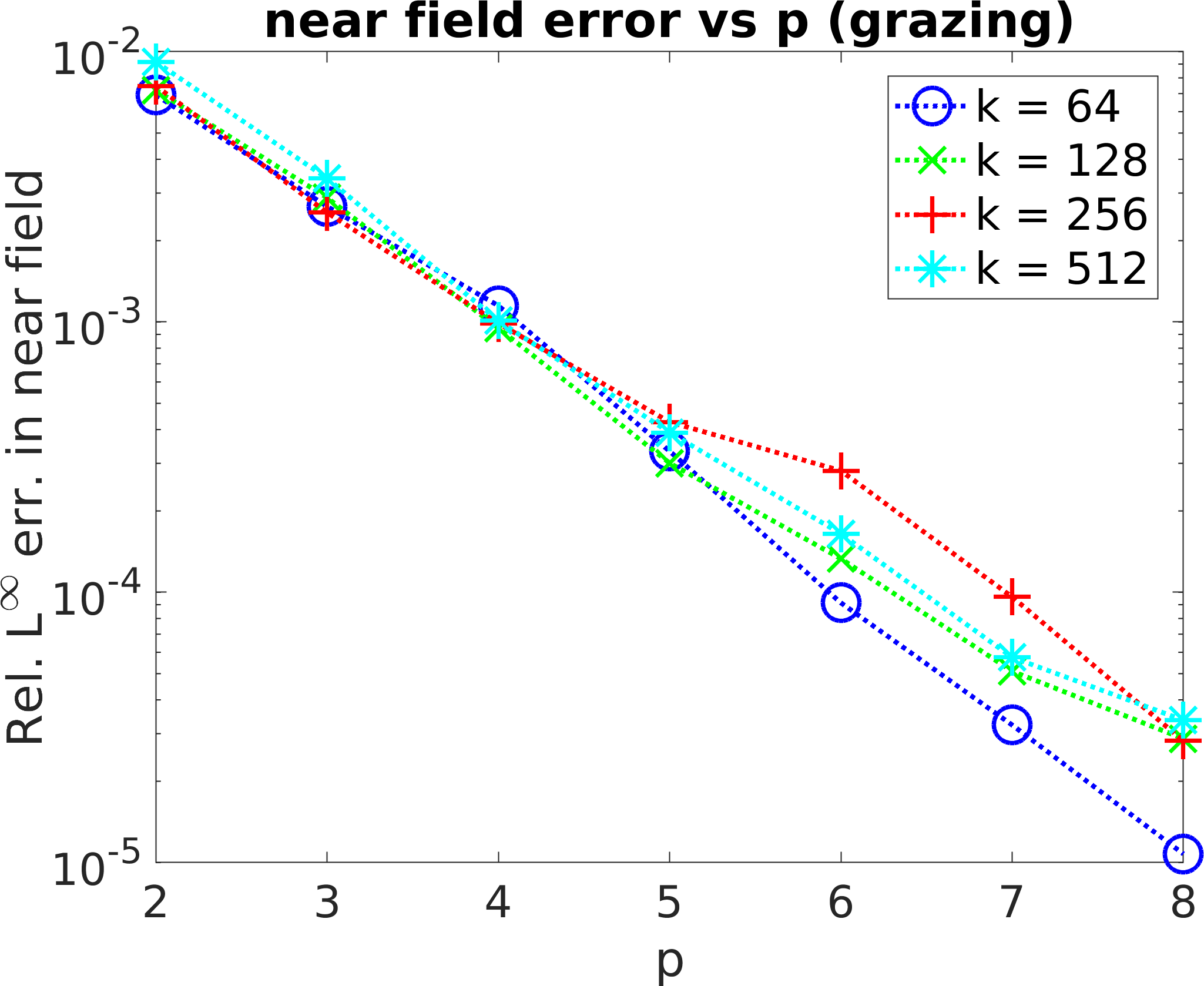}}\hspace{3mm}
	\subfloat[][]{\includegraphics[width=.48\linewidth]{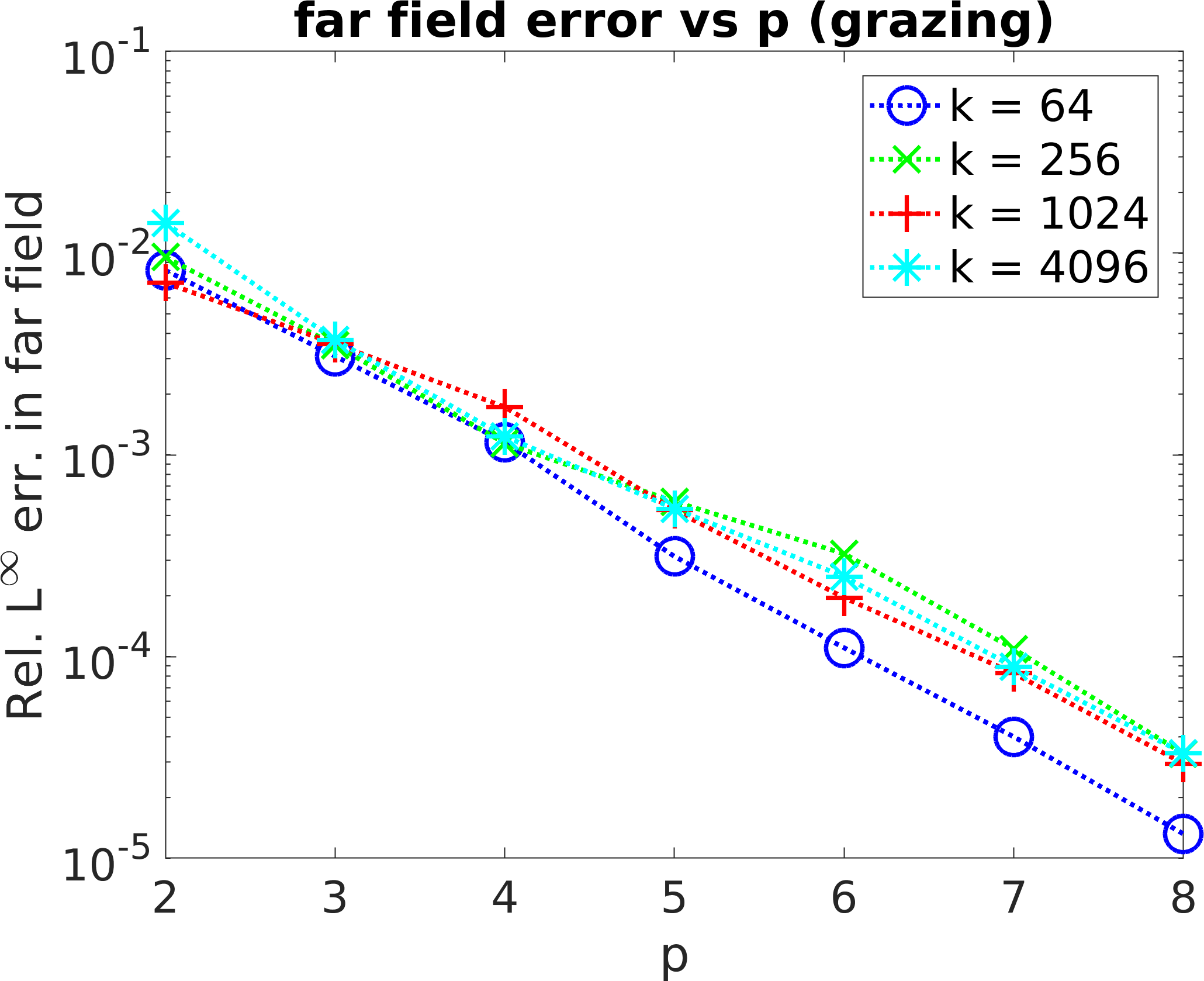}}
	\caption{Analogue of Figs.~\ref{fig:referenceSoln}(a), \ref{fig:performance}(c) and \ref{fig:performance}(d) for the case of grazing incidence.}
	\label{fig:grazing}
\end{figure}

In Fig.~\ref{fig:grazing} we show analogous results for the case of grazing incidence, with $\bd=(0,-1)$. A plot of the corresponding total field is shown in Fig.~\ref{fig:grazing}(a) and plots of the near- and far-field errors are shown in Figs.~\ref{fig:grazing}(b) and \ref{fig:grazing}(c). 
As before we see exponential convergence with increasing $p$, and no significant increase in error for increasing $k$. 
Relative errors are an order of magnitude larger than in the non-grazing example, because in this case there is no geometrical optics component ($\Psi=0$), so the scattered fields are purely diffractive and hence smaller in maximum magnitude.

\resub{Comparing Fig.~\ref{fig:performance}(d) and Fig.~\ref{fig:grazing}(c), we observe that in the far field the accuracy noticeably improves as $k$ increases in the case of non-grazing incidence, but not in the case of grazing incidence. This is due to the fact that the far-field pattern $u^\infty$ has peaks of magnitude proportional to $|d_2|k$ at the values of $\theta$ corresponding to the reflected and shadow directions (cf.\ Fig.~\ref{fig:FFCantor}(a)), associated with the contribution of the geometrical optics term $\Psi$.\footnote{\resub{This contribution is found by substituting $\Psi(\bfx(s))=-2ik|d_2|\re^{\ri k d_1s}$ for $v(s)$ in \eqref{def:FarField}. For example, for $\tGamma=(0,1)$ we compute 
			\[2\ri k |d_2|\int_0^1 \re^{\ri ks(d_1-\cos\theta)}\,\rd s = 2\ri k|d_2|\frac{\re^{\ri k(d_1-\cos\theta)}-1}{\ri k(d_1-\cos\theta)},\]
			which takes the value $2\ri k|d_2|$ when $\cos\theta=d_1$.}} Hence the denominator in the relative $L^\infty$ error will be $O(k)$ for non-grazing incidence, while the numerator is not expected to be comparably large because the contribution from $\Psi$ is computed exactly (up to quadrature error), and hence does not contribute to the discretization error.}

\subsection{Collocation parameters}\label{sec:tuningParams}

Our choices of oversampling parameter $C_\mathrm{OS}=1.25$ and SVD truncation parameter $\epsilon=10^{-8}$ in the previous section were based on the results of extensive numerical testing. In this section we report a small subset of these results, to illustrate the effect of changing these parameters. All results in this section relate to the case $\Gamma=(0,1)\times\{0\}$ with $\bd=(1,-1)/\sqrt{2}$, as in Figs.~\ref{fig:referenceSoln} and \ref{fig:performance}.

\begin{figure}[tp!]
	\subfloat[][]{\includegraphics[width=.47\linewidth]{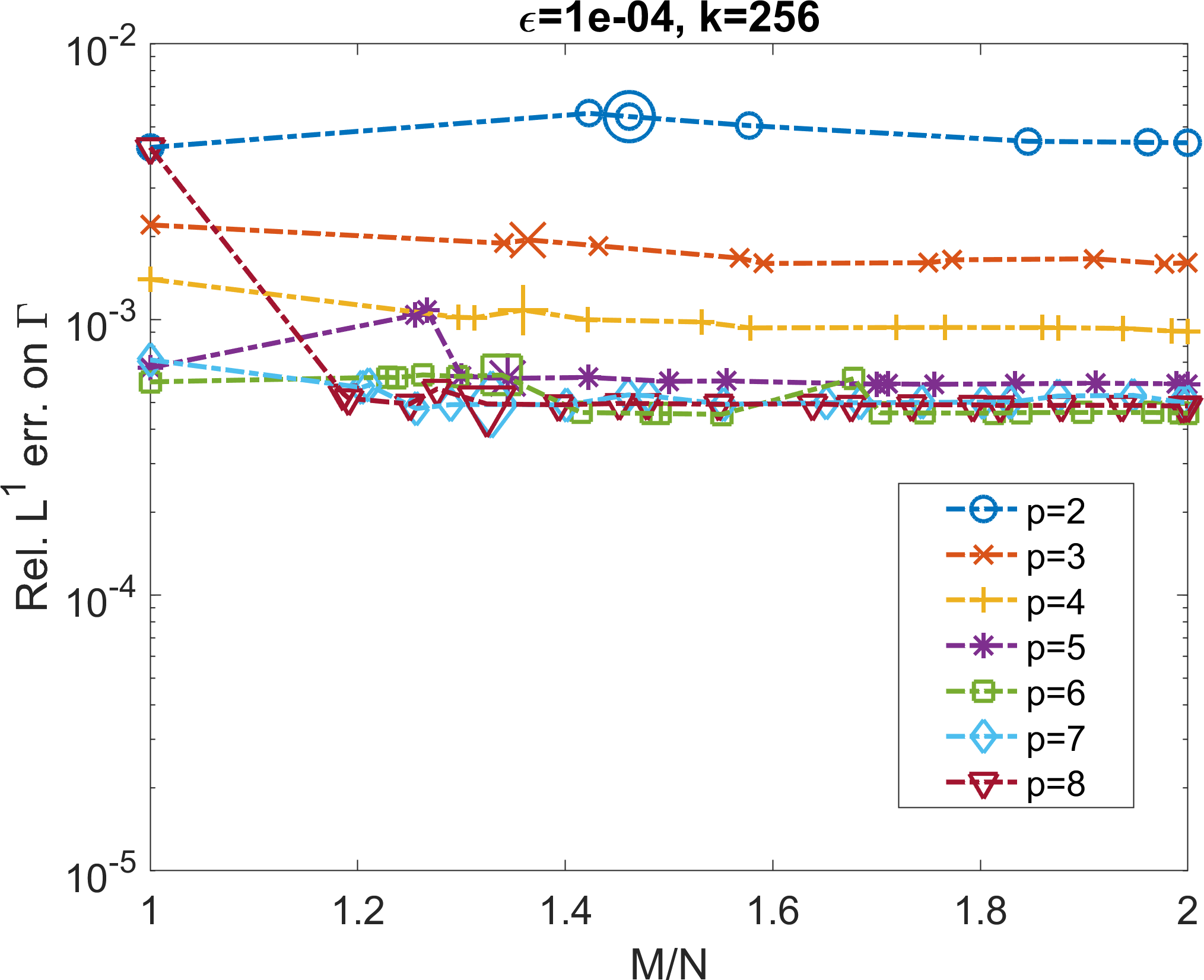}}
	\hspace{3mm}
	\subfloat[][]{\includegraphics[width=.47\linewidth]{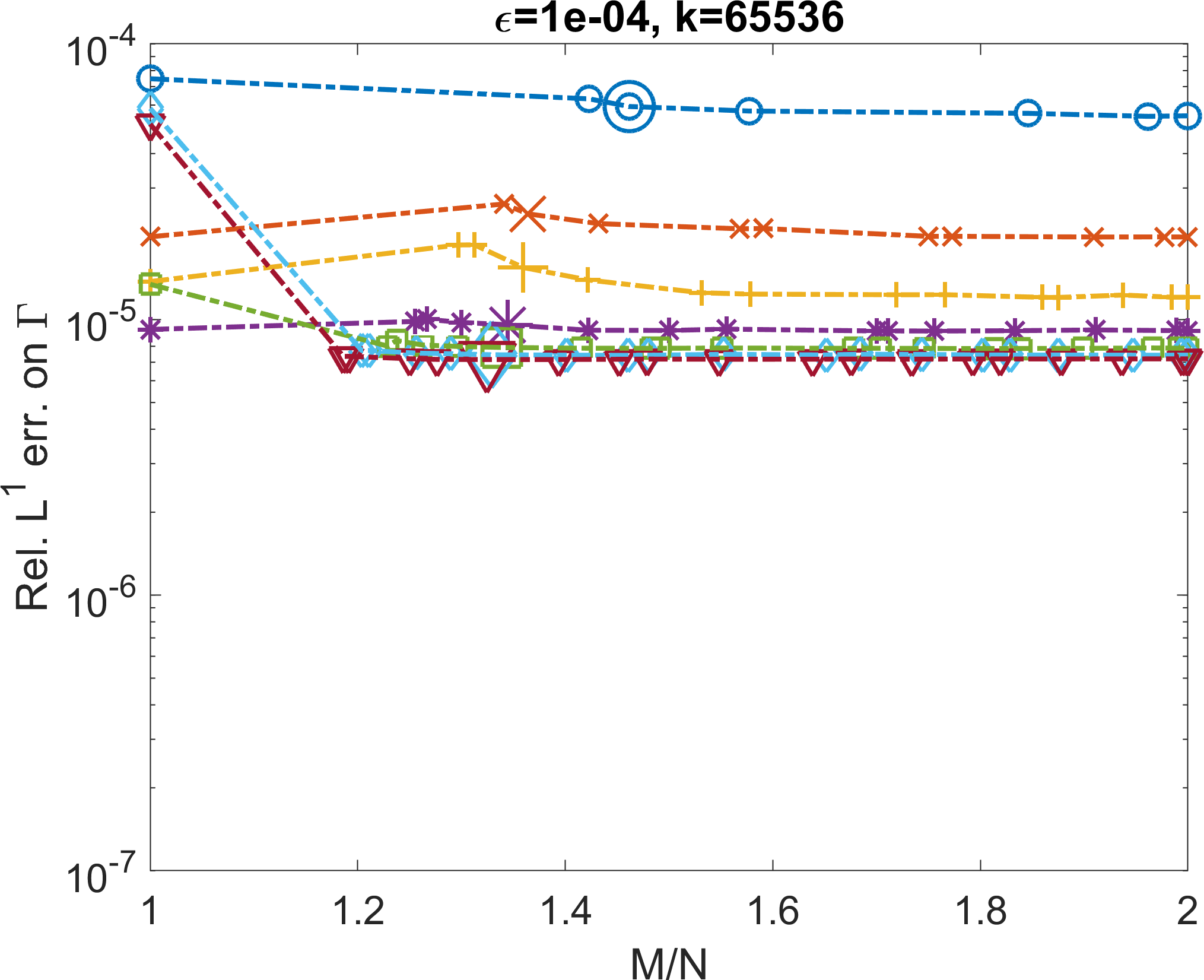}}\\
	\subfloat[][]{\includegraphics[width=.47\linewidth]{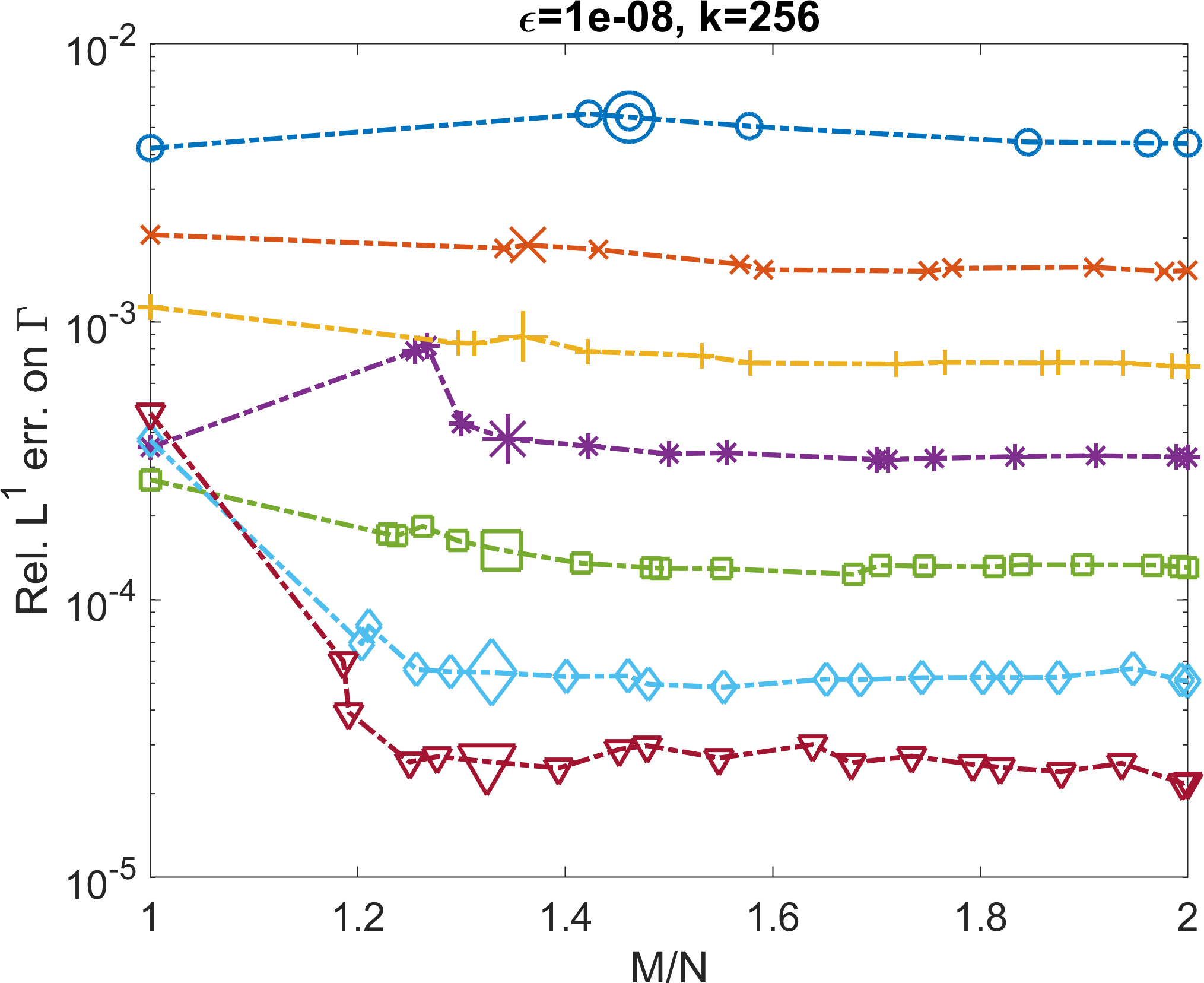}}
	\hspace{3mm}
	\subfloat[][]{\includegraphics[width=.47\linewidth]{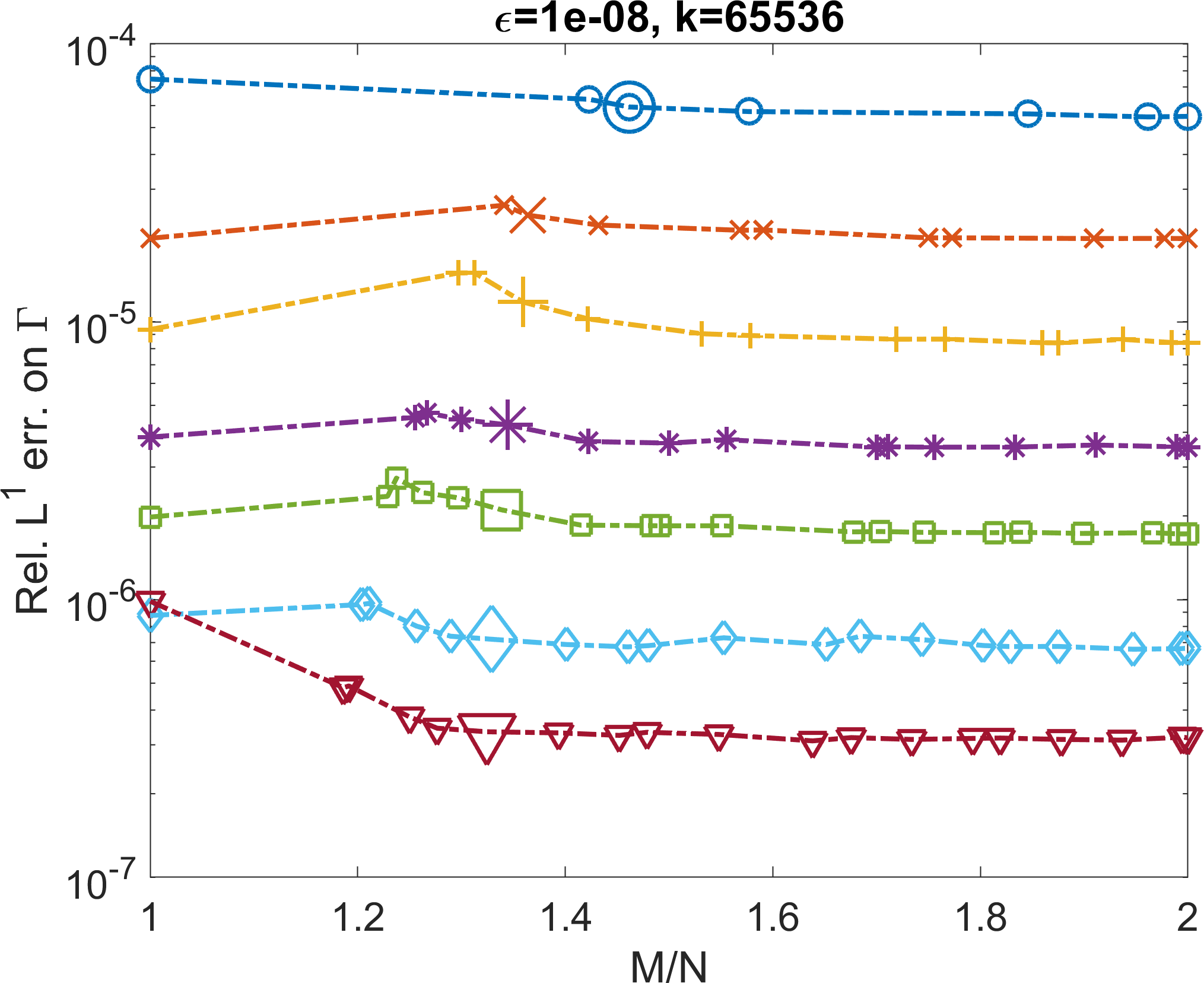}}\\
	\subfloat[][]{\includegraphics[width=.47\linewidth]{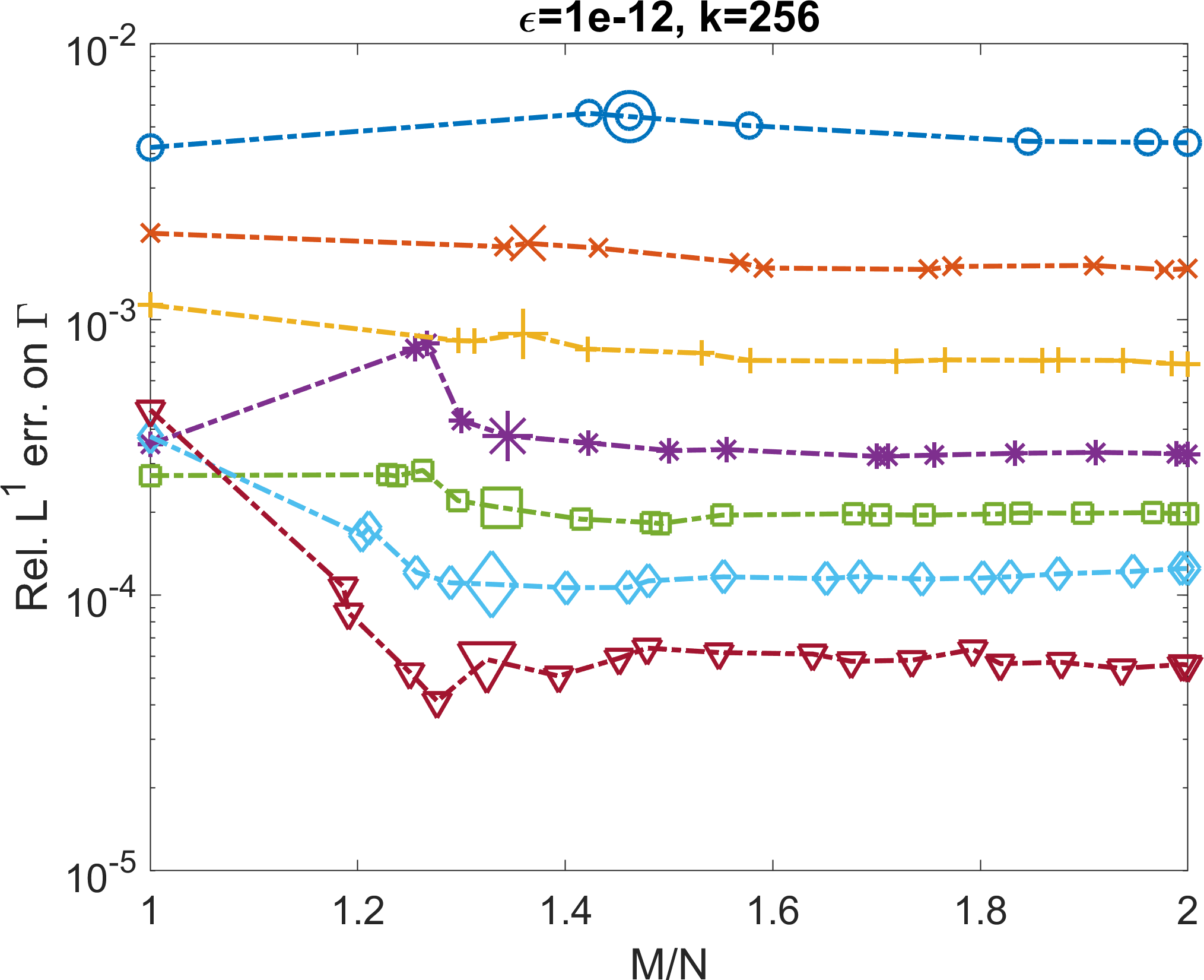}}
	\hspace{3mm}
	\subfloat[][]{\includegraphics[width=.47\linewidth]{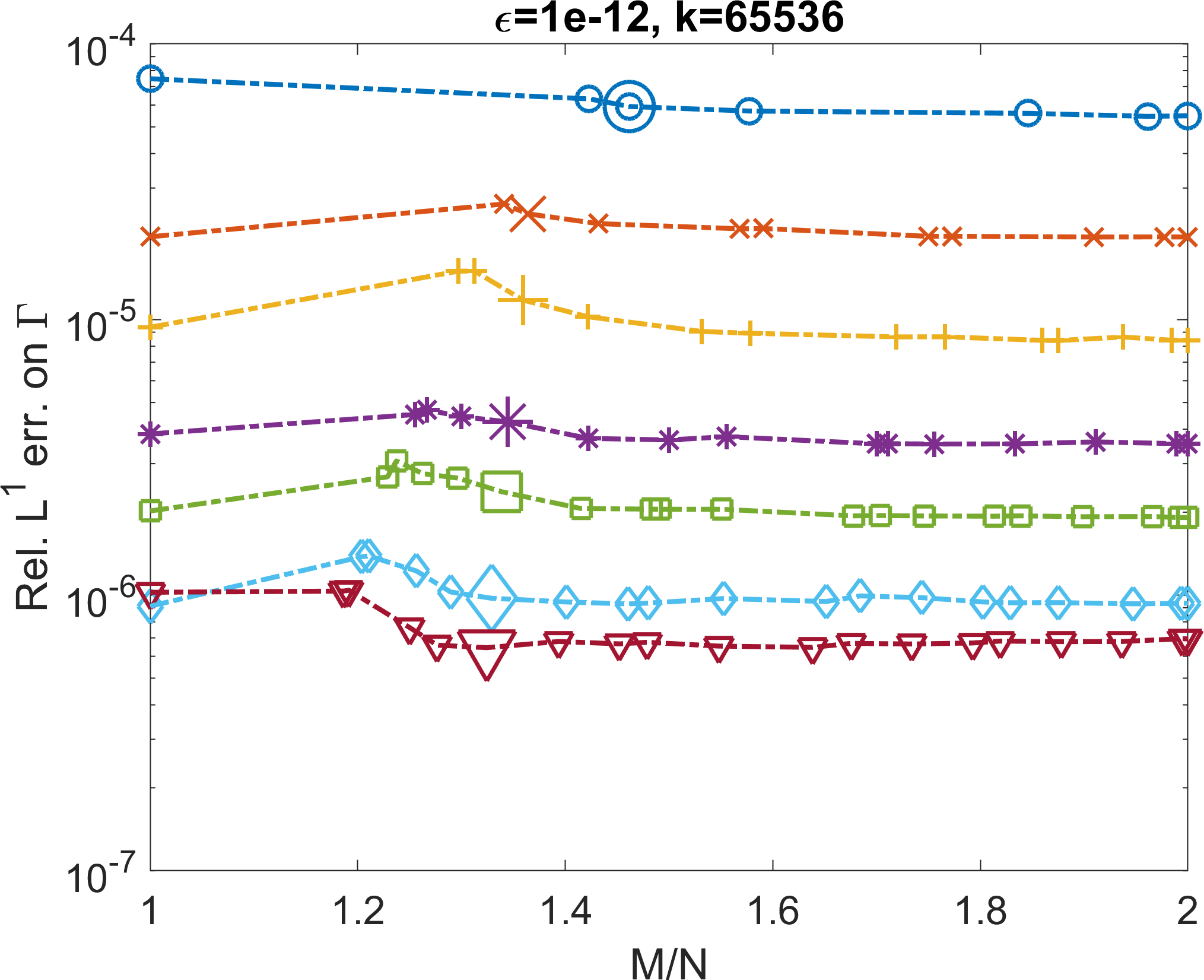}}
	\vspace{-2mm}
	\caption{
		Dependence of solution accuracy on the oversampling parameter $C_\text{OS}$ and the truncation parameter $\epsilon$ for moderate ($k=256$, left-hand plots) and large ($k=65536$, right-hand plots) wavenumbers and varying polynomial degree $p$. The values of $C_\text{OS}$ used are $\{1,1.05,1.1,\ldots,2\}$, but the resulting ratio $M/N$ between the number of collocation points and the number of unknowns is in general slightly larger than $C_\text{OS}$ (recall \rf{eq:MoverN}). 
		The data points corresponding to $C_\mathrm{OS}=1.25$ are shown with large markers.
	}\label{fig:COSeps_compare}
\end{figure}

In Fig.~\ref{fig:COSeps_compare} we present plots showing the dependence of the $L^1(\Gamma)$ solution error, for $p\in\{2,\ldots,8\}$, on the oversampling ratio $M/N$, which is a function of the oversampling parameter $C_\mathrm{OS}$. The left-hand plots correspond to wavenumber $k=256$ and the right-hand plots to $k=65536$, while the top, middle and bottom plots correspond to the SVD truncation thresholds $\epsilon=10^{-4}$, $\epsilon=10^{-8}$ and $\epsilon=10^{-12}$, respectively. In each case tests were run for $C_\mathrm{OS}\in\{1,1.05,1.1,\ldots,2\}$; note that for the non-integer values of $C_\mathrm{OS}$ the resulting values of $M/N$ are in general larger than $C_\mathrm{OS}$ (recall \rf{eq:MoverN}). 
As the reference solution we use the solution with $p=12$ and $C_\mathrm{OS}=2$.

In interpreting these plots we focus first on the results for $C_\mathrm{OS}=1$ ($M/N=1$), which corresponds to a square system (no oversampling). 
The need for oversampling is clear from the fact that in all six plots the errors for $C_\mathrm{OS}=M/N=1$ are not monotonically decreasing with increasing $p$. That is, increasing the dimension of the approximation space $V_N$ does not necessarily lead to a more accurate solution. (Further results for this square case 
can be found in \cite{emile}.) 

By contrast, at least in plots (c)-(f) of Fig.~\ref{fig:COSeps_compare}, for all fixed values of $C_\mathrm{OS}>1$ considered, the error is monotonically decreasing with increasing $p$. Since computational cost is proportional to $M$, it is desirable to take $C_\mathrm{OS}$ as small as possible while maintaining this monotonicity. But for $C_\mathrm{OS}$ close to $1$ there are still visible instabilities. On the basis of our experiments it seems that $C_\mathrm{OS}\geq 1.25$ is sufficient to ensure stability for all the problems we considered. (The data points corresponding to the choice $C_\mathrm{OS}=1.25$ are shown with larger markers to make this transition more clear.) Significantly, and perhaps surprisingly, the amount of oversampling required to stabilise the method appears to be independent of the wavenumber $k$, as one sees, for example, by comparing Figs.~\ref{fig:COSeps_compare}(c) (moderate $k$) and \ref{fig:COSeps_compare}(d) (larger $k$), and recalling the convergence plots in Figs.~\ref{fig:performance}(a) and \ref{fig:performance}(b), which reach wavenumber $k=262,144$ with no discernable loss of stability. 

To interpret the dependence of the results on the SVD truncation parameter $\epsilon$, it is instructive to recall Lemma \ref{lem:svd}, which, while not directly providing information about the error in our boundary solution (as explained in \S\ref{sec:SVD}), still provides some intuition. On the one hand, since the argument of the infimum on the right-hand side of \rf{eq:svd} is a linearly increasing function of $\epsilon$, we expect that if $\epsilon$ is taken to be too large, the discrete residual for the SVD solution may be %unacceptably 
large, leading to a loss in accuracy of the boundary solution. This effect can be clearly observed by comparing Figs.~\ref{fig:COSeps_compare}(c) and \ref{fig:COSeps_compare}(d) ($\epsilon=10^{-8}$), where we see clear exponential convergence as $p$ increases, with Figs.~\ref{fig:COSeps_compare}(a) and \ref{fig:COSeps_compare}(b) ($\epsilon=10^{-4}$), where errors appear to stagnate around $p=5$.
On the other hand, looking again at \rf{eq:svd}, we expect that if $\epsilon$ is taken to be too small, the SVD solver may select a solution with a large coefficient norm, in which case our numerical solution may suffer from spurious oscillations between the collocation points, or other numerical instabilities that might increase the $L^1(\Gamma)$ error. This degradation in accuracy for small $\epsilon$ can be observed by comparing Figs.~\ref{fig:COSeps_compare}(c) and \ref{fig:COSeps_compare}(d) ($\epsilon=10^{-8}$), with Figs.~\ref{fig:COSeps_compare}(e) and \ref{fig:COSeps_compare}(f) ($\epsilon=10^{-12}$). 
Based on extensive testing we found that the choice $\epsilon=10^{-8}$ gave satisfactory performance for all the examples we considered. 

\subsection{Scattering by a Cantor set}\label{sec:Cantor}

In this section we present an application of our HNA BEM to high frequency scattering by 
the 
middle-third 
Cantor set. 
The mathematical analysis of acoustic scattering by fractal screens, of which the Cantor set is an example, 
was developed recently in \cite{ChHe:18} and \cite{ChHeMoBe:19}, with the latter providing a rigorous convergence theory for an approximation strategy based on conventional Galerkin BEM using piecewise constant basis functions, along with numerical results for low frequency scattering in both two and three dimensions by a range of fractal screens. 
For two-dimensional scattering by the Cantor set, the HNA BEM presented in the current paper allows us to investigate the high frequency regime, and to 
our knowledge this is the first time that high frequency BEM simulations of scattering by fractal screens have been presented in the literature. 
Fractals, 
which one might define loosely as ``sets possessing geometric detail on every lengthscale'', 
arise in numerous scattering applications, including the human lung and other dendritic structures in medical science \cite{JoSe:11}, 
snowflakes, ice crystals and other atmospheric particles in climate science \cite{StWe:15}, 
fractal antennas for electromagnetic wave transmission/reception \cite{GhSiKa:14}, fractal piezoelectric ultrasound transducers \cite{MuWa:11}, and fractal aperture problems in laser physics \cite{Ch:16}; for further references see \cite{ChHeMoBe:19}.
The case of high frequency scattering by fractals is particularly intriguing, because even though the diameter of the scatterer may be large compared to the wavelength, the scatterer always possesses sub-wavelength structure, so standard ray-based high frequency approximations do not apply. 
In his 1979 paper on random phase screens \cite{Berry79}, 
Berry describes this case as ``a new regime in wave physics" in which the scattered waves adopt ``unfamiliar forms that should be studied in their own right".% \cite{Berry79}. 

In practical applications and numerical simulations one generally works with ``prefractal'' approximations in which the fractal structure is truncated at some level. 
The standard prefractals for the middle-third Cantor set are defined as follows. 
Letting $\Gamma^{(0)}=(0,1)\times\{0\}$, for $j\in \N$ we construct the $j$th 
prefractal $\Gamma^{(j)}$ by removing the (closed) middle third from each disjoint interval making up the previous prefractal $\Gamma^{(j-1)}$, so that e.g.\ $\Gamma^{(1)}=((0,1/3)\cup (2/3,1))\times\{0\}$. 
That the BIE solutions for sound-soft scattering on the prefractals converge to a well-defined non-zero limiting BIE solution 
on the underlying fractal 
was proved in \cite[Example 8.2]{ChHe:18} (and see also \cite[Prop.~6.1]{ChHeMoBe:19}).

\begin{figure}
	\subfloat[][]{\includegraphics[width=.4\linewidth]{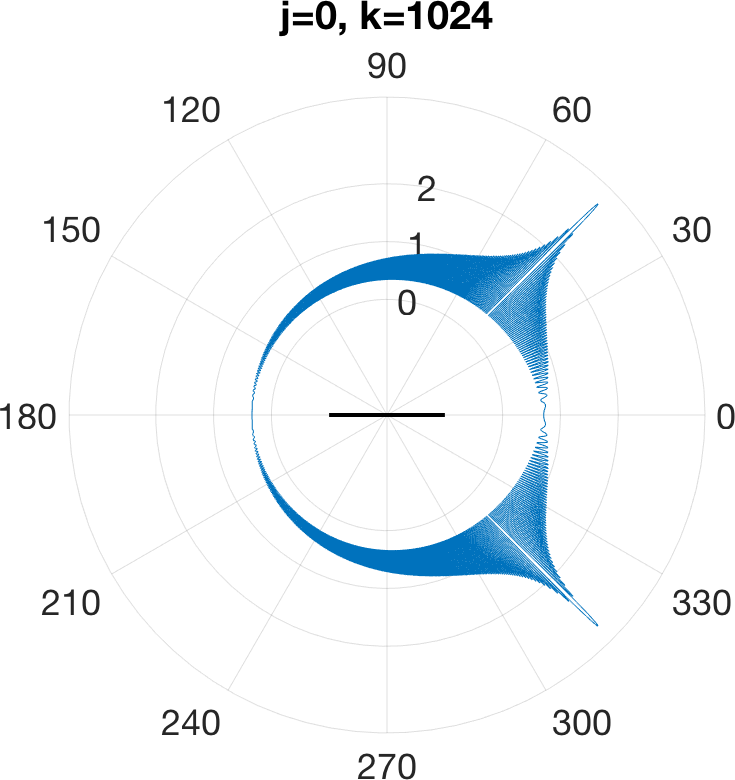}} \hspace{12mm}
	\subfloat[][]{\includegraphics[width=.4\linewidth]{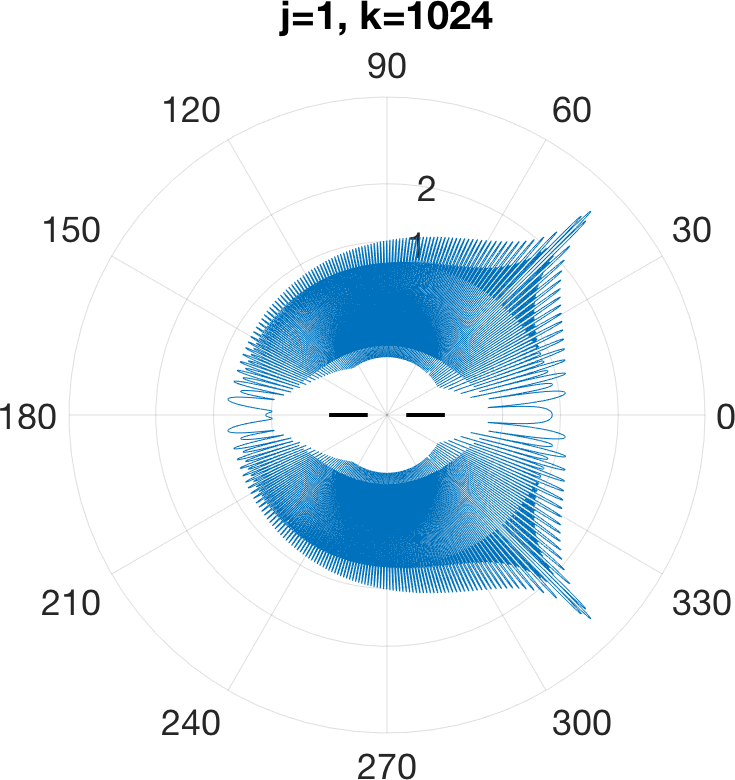}}\\
	\subfloat[][]{\includegraphics[width=.4\linewidth]{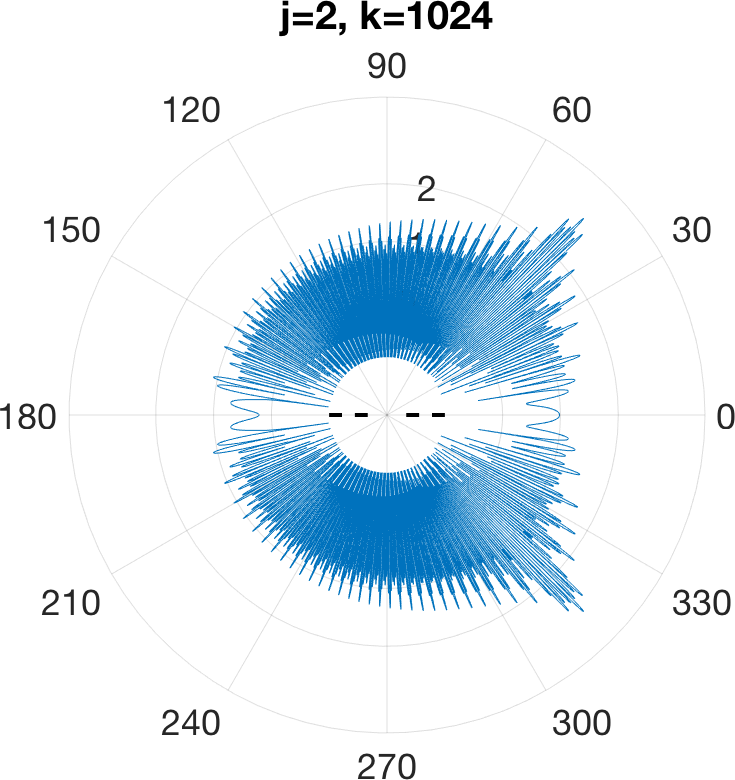}}\hspace{12mm}
	\subfloat[][]{\includegraphics[width=.4\linewidth]{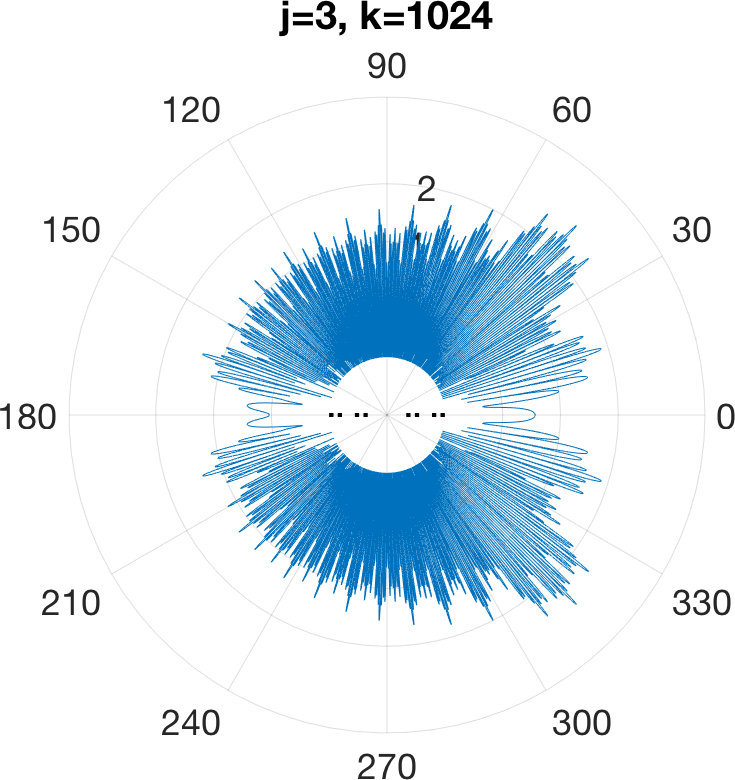}}\\
	\subfloat[][]{\includegraphics[width=.4\linewidth]{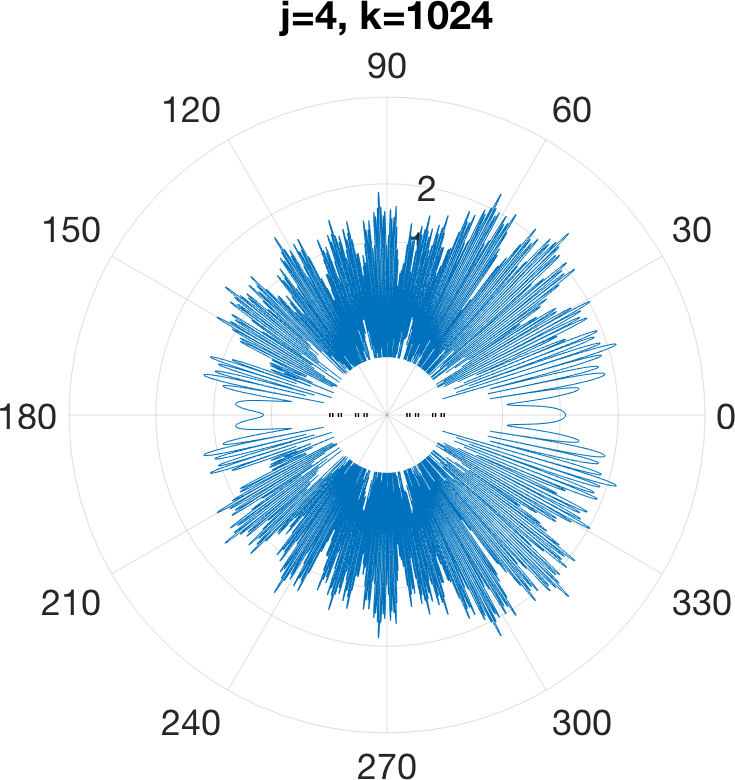}}\hspace{12mm}
	\subfloat[][]{\includegraphics[width=.4\linewidth]{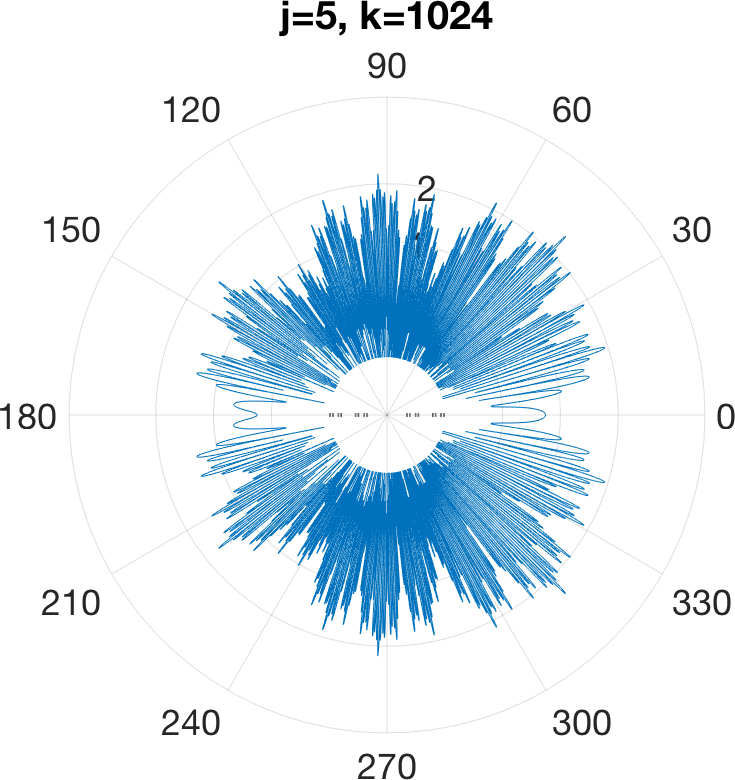}}
	\caption{Radial plots of $\log_{10}|u^\infty_8(\theta)|$ for $k=1024$ and prefractal levels  $j=0,\ldots ,5$, for scattering by the middle-third Cantor set.}\label{fig:FFCantor}
\end{figure}

\begin{figure}
	\centering \includegraphics[width=.5\linewidth]{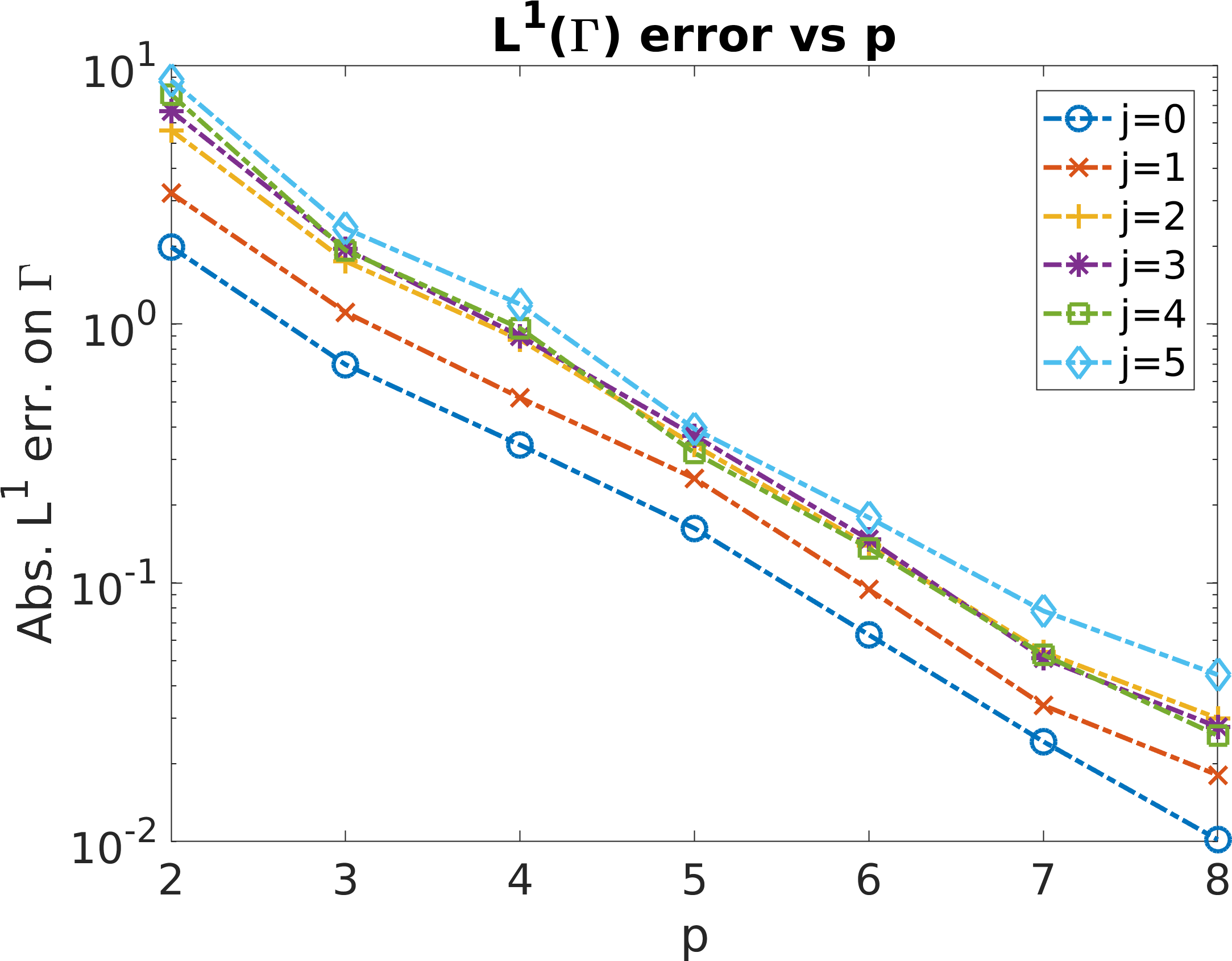}
	\caption{Convergence of the boundary solution with increasing $p$ for $k=1024$ and prefractal levels $j=0,\ldots ,5$, for scattering by the middle-third Cantor set.}
	\label{fig:CantorErr}
\end{figure}

A plot of the total field for scattering by $\Gamma^{(2)}$ with incident direction $\bd=(1,-1)/\sqrt{2}$, computed using our collocation HNA BEM, was presented in Fig.~\ref{fig:mutliscreenplot}. 
In Fig.~\ref{fig:FFCantor} we present radial plots of \[\max\{\log_{10}|u^\infty_8|,-1\},\] for the first six prefractals $\Gamma^{(0)}, \ldots \Gamma^{(5)}$, with $k=1024$ and $\bd=(1,-1)/\sqrt{2}$. 
The log scales are truncated at $-1$ to leave space in the centre of the plots for the corresponding prefractals to be plotted. 
Corresponding $L^1(\Gamma)$ convergence plots are presented in Fig.~\ref{fig:CantorErr}, where, for each $j=0, \ldots, 5$ we take as reference solution the corresponding solution with $p=10$. 
The fact that the absolute $L^1(\Gamma)$ error in our boundary solution with $p=8$ is smaller than $10^{-1}$ for all $j=0, \ldots, 5$ implies that our far-field plots in Fig.~\ref{fig:FFCantor} are also accurate to $10^{-1}$ absolute error, since
\begin{align}
	\label{eqn:BdytoFF}
	|u^\infty_p|\leq \|v_p\|_{L^1(\Gamma)}.
\end{align}
We note that at $k=1024$ the 
lowest order prefractal $\Gamma^{(0)}$ is approximately 160 wavelengths long, and each component of the highest order prefractal $\Gamma^{(5)}$ is approximately one wavelength long. 
While the far-field patterns for $\Gamma^{(4)}$ and $\Gamma^{(5)}$ are similar, they are clearly not identical, which is to be expected since in refining $\Gamma^{(4)}$ to $\Gamma^{(5)}$ we are making wavelength-scale changes to the scatterer geometry. 
Hence for $\Gamma^{(5)}$ we are still in the ``pre-asymptotic'' regime with regard to convergence to the solution for scattering by the limiting fractal screen. 
(The asymptotic regime is considered in \cite{ChHeMoBe:19} using a conventional BEM approach.) 
We emphasize that our HNA method can handle much larger wavenumbers than $k=1024$, but at very large wavenumbers the far-field plots are so oscillatory they are difficult to interpret, so we do not present them here. 

\begin{acknowledgements}
	The authors acknowledge support from KU Leuven project C14/15/055 (A.~Gibbs) and EPSRC project EP/S01375X/1 (D.~Hewett and A.~Gibbs).
\end{acknowledgements}

\bibliographystyle{spmpsci}
\bibliography{HNA_bib}

\end{document}